\documentclass[11pt]{article}

\usepackage{cite}             
\usepackage{latexsym}         
\usepackage{amsmath}          
\usepackage{amssymb}          
\usepackage{amsxtra}          
\usepackage[latin1]{inputenc} 
\usepackage{eucal}            
\usepackage{longtable}        
\usepackage{exscale}          
\usepackage{bbm}              
\usepackage{theorem}          
\usepackage{paralist}         
\usepackage[curve]{xypic}
\usepackage{diagxy}           
\usepackage{stmaryrd}         

\title{Classification of Invariant Star Products up to Equivariant Morita
Equivalence on Symplectic Manifolds}

\author{\textbf{Stefan Jansen},
  \textbf{Nikolai Neumaier},
  \textbf{Gregor Schaumann},
  \textbf{Stefan Waldmann}\thanks{Corresponding author:
    Stefan.Waldmann@physik.uni-freiburg.de}
  \\[0.1cm]
  Fakult{\"a}t f{\"u}r Mathematik und Physik\\
  Albert-Ludwigs-Universit{\"a}t Freiburg\\
  Physikalisches Institut\\
  Hermann-Herder-Stra{\ss}e 3\\
  D 79104 Freiburg\\
  Germany}

\date{April 2010}

\renewcommand{\mathbb}[1]{\mathbbm{#1}} 

\newcommand{\cc}[1]      {\overline{{#1}}}              
\newcommand{\id}         {\operatorname{\mathsf{id}}}   
\newcommand{\image}      {\operatorname{{\mathrm{im}}}} 
\newcommand{\Lie}        {\operatorname{\mathcal{L}}}    
\newcommand{\ad}         {\operatorname{\mathrm{ad}}}    
\newcommand{\Hom}        {\operatorname{\mathsf{Hom}}}   
\newcommand{\End}        {\operatorname{\mathsf{End}}}   
\newcommand{\SP}[1]      {\left\langle{#1}\right\rangle} 
\newcommand{\ring}[1]    {\mathsf{#1}}                 
\newcommand{\Unit}       {\mathbb{1}}                  
 
\newcommand{\cl}         {\mathrm{cl}}                    
\newcommand{\I}          {\mathrm{i}}

\newcommand{\D}          {\operatorname{\mathrm{d}}} 
\newcommand{\lie}[1]     {\mathfrak{#1}}
\newcommand{\acts}       {\mathbin{\triangleright}}
\newcommand{\actsp}      {\mathbin{\triangleright'}}
\newcommand{\ccacts}     {\mathbin{\cc{\triangleright}}}
\newcommand{\sweedler}[1] {{\scriptscriptstyle{(#1)}}}

\newcommand{\Sym}        {\mathrm{S}}
\newcommand{\Schnitte}   {\Gamma^\infty}
\newcommand{\Formen}     {\Omega}
\newcommand{\ins}        {\operatorname{\mathrm{i}}}
\newcommand{\HdR}        {\mathrm{H}_{\scriptscriptstyle\mathrm{dR}}}
\newcommand{\deform}[1]    {\boldsymbol{#1}}
\newcommand{\Fkt}          {C^{\infty}(M)}
\newcommand{\Bimod}[5] {\sideset{^{\scriptscriptstyle{#1}}_{\scriptscriptstyle{#2}}}{^{\scriptscriptstyle{#4}}_{\scriptscriptstyle{#5}}}{\operatorname{#3}}}

\newcommand{\EA}   {\Bimod{}{}{\mathcal{E}}{}{\mathcal{A}}}
\newcommand{\BEA}  {\Bimod{}{\mathcal{B}}{\mathcal{E}}{}{\mathcal{A}}}
\newcommand{\BEpA} {\Bimod{}{\mathcal{B}}{\mathcal{E}}{\prime}{\mathcal{A}}}
\newcommand{\CFB}  {\Bimod{}{\mathcal{C}}{\mathcal{F}}{}{\mathcal{B}}}
\newcommand{\AccEB}{\Bimod{}{\mathcal{A}}{\cc{\mathcal{E}}}{}{\mathcal{B}}}
\newcommand{\AccEA}{\Bimod{}{\mathcal{A}}{\cc{\mathcal{E}}}{}{\mathcal{A}}}

\newcommand{\AtildeEB}{\Bimod{}{\mathcal{A}}{\widetilde{\mathcal{E}}}{}{\mathcal{B}}}
\newcommand{\AAA}  {\Bimod{}{\mathcal{A}}{\mathcal{A}}{}{\mathcal{A}}}
\newcommand{\BBB}  {\Bimod{}{\mathcal{B}}{\mathcal{B}}{}{\mathcal{B}}}
\newcommand{\AEA}  {\Bimod{}{\mathcal{A}}{\mathcal{E}}{}{\mathcal{A}}}

\newcommand{\IP}[4]{{\,}_{\scriptscriptstyle{#2}\!\!}\left\langle{{#1}}\right\rangle^{\scriptscriptstyle{#3}}_{\scriptscriptstyle{#4}}}

\newcommand{\SPA}[1]     {\IP{{#1}}{}{}{\mathcal{A}}}
\newcommand{\SPEA}[1]    {\IP{{#1}}{}{\mathcal{E}}{\mathcal{A}}}

\newcommand{\SPFB}[1]    {\IP{{#1}}{}{\mathcal{F}}{\mathcal{B}}}
\newcommand{\SPFEA}[1]   {\IP{{#1}}{}{\mathcal{F}\otimes\mathcal{E}}{\mathcal{A}}}
\newcommand{\BSP}[1]     {\IP{{#1}}{\mathcal{B}}{}{}}

\newcommand{\tensor}[1][{}]{\mathbin{\otimes_{\scriptscriptstyle{#1}}}}
\newcommand{\tensM}[1][{}] {\mathbin{\widehat{\otimes}_{\scriptscriptstyle{#1}}}}

\newcommand{\BiMod}       {\operatorname{\mathsf{Bimod}}}
\newcommand{\starBiMod}   {\sideset{}{^{*}}{\operatorname{\mathsf{Bimod}}}}
\newcommand{\StrBiMod}    {\sideset{}{^{\mathrm{str}}}{\operatorname{\mathsf{Bimod}}}}
\newcommand{\BiModH}       {\operatorname{\mathsf{Bimod}_H}}
\newcommand{\starBiModH}   {\sideset{}{^{*}_H}{\operatorname{\mathsf{Bimod}}}}
\newcommand{\StrBiModH}    {\sideset{}{^{\mathrm{str}}_H}{\operatorname{\mathsf{Bimod}}}}

\newcommand{\Pic}      {\operatorname{\mathsf{Pic}}}
\newcommand{\PicH}     {\sideset{}{_{H}}{\operatorname{\mathsf{Pic}}}}
\newcommand{\StrPic}   {\sideset{}{^{\mathrm{str}}}{\operatorname{\mathsf{Pic}}}}
\newcommand{\starPic}  {\sideset{}{^*}{\operatorname{\mathsf{Pic}}}}
\newcommand{\StrPicH}  {\sideset{}{^{\mathrm{str}}_{H}}{\operatorname{\mathsf{Pic}}}}
\newcommand{\starPicH} {\sideset{}{^*_{H}}{\operatorname{\mathsf{Pic}}}}

\newcommand{\starAut}  {\sideset{}{^*}{\operatorname{\mathsf{Aut}}}}
\newcommand{\starIsoH} {\sideset{}{^*_{H}}{\operatorname{\mathsf{Iso}}}}
\newcommand{\starAutH} {\sideset{}{^*_{H}}{\operatorname{\mathsf{Aut}}}}
\newcommand{\AutH}     {\sideset{}{_{H}}{\operatorname{\mathsf{Aut}}}}

\newcommand{\SPicH}      {\sideset{}{_{H}}{\operatorname{\mathsf{SPic}}}}

\newcommand{\StrSPicH}   {\sideset{}{^{\mathrm{str}}_H}{\operatorname{\mathsf{SPic}}}}

\newcommand{\starSPicH}  {\sideset{}{^*_{H}}{\operatorname{\mathsf{SPic}}}}

\newcommand{\GR}[2]             {\ensuremath{\mathsf{#1}(#2)}}
\newcommand{\GRn}[2]            {\ensuremath{\mathsf{#1}_{0}(#2)}} 
\newcommand{\Zentrum}[2][{}]  {\sideset{}{_{\scriptscriptstyle {#1}}}{\operatorname{\mathfrak{Z}}}(#2)}

\newcommand{\Glha}       {\GR{Gl}{H, \alg{A}}}
\newcommand{\Glhb}       {\GR{Gl}{H, \alg{B}}}
\newcommand{\Uha}        {\GR{U}{H,\alg{A}}}
\newcommand{\Uhb}         {\GR{U}{H,\alg{B}}}
\newcommand{\Gloha}       {\GRn{Gl}{H, \alg{A}}}
\newcommand{\Glohb}       {\GRn{Gl}{H, \alg{B}}}
\newcommand{\Uoha}        {\GRn{U}{H,\alg{A}}}
\newcommand{\Uohb}         {\GRn{U}{H,\alg{B}}}
\newcommand{\Glhza}       {\GR{Gl}{H, \Zentrum{\alg{A}}}}
\newcommand{\Uhza}        {\GR{U}{H,\Zentrum{\alg{A}}}}

\newcommand{\Glza}       {\GR{GL}{\Zentrum{\alg{A}}}}
\newcommand{\Uza}        {\GR{U}{\Zentrum{\alg{A}}{}}}
\newcommand{\Glzb}      {\GR{GL}{\Zentrum{\alg{B}}{}}}
\newcommand{\Uzb}        {\GR{U}{\Zentrum{\alg{B}}{}}}

\newcommand{\alg}[1]       {\mathcal{#1}}

\newcommand{\starMM}       {\operatorname{\mathrm{MoMa}}^*}
\newcommand{\MM}           {\operatorname{\mathrm{MoMa}}}

\theoremheaderfont{\normalfont\bfseries}
\theorembodyfont{\itshape}
\newtheorem{lemma}{Lemma}[section]
\newtheorem{proposition}[lemma]{Proposition}
\newtheorem{theorem}[lemma]{Theorem}

\newtheorem{definition}[lemma]{Definition}

\theorembodyfont{\rmfamily}

\newtheorem{example}[lemma]{Example}
\newtheorem{remark}[lemma]{Remark}

\makeatletter
\newcommand\qedsymbol{\hbox{$\boxempty$}}
\newcommand\qed{\relax\ifmmode\boxempty\else
  {\unskip\nobreak\hfil\penalty50\hskip1em\null\nobreak\hfil\qedsymbol
  \parfillskip=\z@\finalhyphendemerits=0\endgraf}\fi}
\makeatother
\newenvironment{proof}[1][{}]{\par\noindent Proof{#1}. }{\qed}

\newenvironment{theoremlist}{\begin{enumerate}}{\end{enumerate}}

\newenvironment{lemmalist}{\begin{enumerate}}{\end{enumerate}}
\newenvironment{propositionlist}{\begin{enumerate}}{\end{enumerate}}
\newenvironment{definitionlist}{\begin{enumerate}}{\end{enumerate}}

\numberwithin{equation}{section}

\setdefaultenum{\itshape i.)}{a.)}{I.)}{A.)}
\newcommand{\refitem}[1] {~\textit{\ref{#1}.)}}

\begin{document}

\maketitle

\begin{abstract}
    In this paper we investigate equivariant Morita theory for
    algebras with momentum maps and compute the equivariant Picard
    groupoid in terms of the Picard groupoid explicitly. We consider
    three types of Morita theory: ring-theoretic equivalence,
    $^*$-equivalence and strong equivalence. Then we apply these
    general considerations to star product algebras over symplectic
    manifolds with a Lie algebra symmetry. We obtain the full
    classification up to equivariant Morita equivalence.
\end{abstract}

\newpage

%
%

\tableofcontents

%
%

\section{Introduction}
\label{sec:Intro}

The aim of this paper is two-fold: on the one hand we continue the
study of equivariant Morita theory of algebras with respect to Hopf
algebra actions as initiated in \cite{jansen.waldmann:2006a}. On the
other hand, we consider the particular case of star product algebras
over symplectic manifolds with symmetries coming from Lie algebra
actions on the underlying manifold.

On the general algebraic side, the algebras in question carry a left
action of a fixed Hopf algebra $H$. We are interested in the groupoid
of all isomorphism classes of Morita equivalence bimodules with
compatible action of $H$. Beside this purely algebraic case, we are
interested in the situation that the algebras are defined over a ring
$\ring{C} = \ring{R}(\I)$ with an ordered ring $\ring{R}$ and $\I^2 =
-1$ such that they carry a $^*$-involution. In this case, also the
Hopf algebra is required to be a Hopf $^*$-algebra and the action is
required to be a $^*$-action. Moreover, for the equivalence bimodules
one demands algebra-valued inner products either with or without
complete positivity. This gives the three flavours of Morita theory we
are interested in: ring-theoretic Morita theory, $^*$-Morita theory,
and strong Morita theory. We base our study of Morita theory on the
corresponding Picard groupoids.  While in \cite{jansen.waldmann:2006a}
the general interplay of $H$-equivariance with these three types of
Morita theory was studied, we focus on a more particular situation
where the action of $H$ allows a \emph{momentum map}, i.e. an algebra
homomorphism $J \longrightarrow \mathcal{A}$ which induces the action
as an inner action. In this case, the in general very difficult
question whether Morita equivalence implies also equivariant Morita
equivalence can be answered positively in a trivial way. Even more
precisely, we obtain that the canonical groupoid morphisms of
forgetting the equivariance
\begin{equation}
    \label{eq:PicHtoPicForgetAllStuff}
    \PicH \longrightarrow \Pic,
    \quad
    \starPicH \longrightarrow \starPic,
    \quad
    \textrm{and}
    \quad
    \StrPicH \longrightarrow \StrPic
\end{equation}
are not only surjective but possess a canonical right inverse
constructed from the momentum maps. The kernels of
\eqref{eq:PicHtoPicForgetAllStuff} for all three flavours of Morita
theory were studied in detail in
\cite{jansen.waldmann:2006a}. Moreover, we show that the existence of
a momentum map is an invariant under $H$-equivariant Morita
equivalence. Thus the restriction of the Picard groupoids to this
particular class of actions and algebras is reasonable as we would not
leave these connected components of the groupoids anyway.

On the star product side we consider $\lie{g}$-invariant star products
on symplectic manifolds with a symplectic action of a
finite-dimensional Lie algebra $\lie{g}$. The action may or may not
arise from a corresponding Lie group action. The only technical
condition we need throughout is that the action preserves a
connection, and consequently  also preserves a symplectic torsion-free
connection. In this case, one has strong existence and classification
results of such $\lie{g}$-invariant star products according to
\cite{bertelson.bieliavsky.gutt:1998a} by a $\lie{g}$-invariant
characteristic class $c^{\lie{g}}(\cdot)$ in the formal series in the
$\lie{g}$-invariant second deRham cohomology of $M$.  A particularly
nice situation is obtained if in addition the classical
$\lie{g}$-action allows for a \emph{classical momentum map}. Then also
the star products which allow a \emph{quantum momentum map} have been
studied and classified according to
\cite{mueller-bahns.neumaier:2004a, gutt.rawnsley:2003a}. We show that
the notion of quantum momentum maps is equivalent to the momentum maps
in the above Hopf algebraic approach to equivariant Morita theory if
the Hopf algebra is chosen to be a ``rescaled'' version of the
universal enveloping algebra of $\lie{g}$.

For general star products on symplectic manifolds the classification
up to Morita equivalence (and also up to $^*$-Morita equivalence and
strong Morita equivalence) is available from
\cite{bursztyn.waldmann:2002a}. Thus the task is to combine both
classification results to get the classification of
$\lie{g}$-invariant star products up to $\lie{g}$-equivariant Morita
equivalence. Here we obtain the main result of this paper:
\begin{theorem}
    \label{theorem:TheRealStuff}%
    Suppose there exists a $\lie{g}$-invariant connection on a
    symplectic manifold $M$ and let $\star$ and $\star'$ be two
    $\lie{g}$-invariant (Hermitian) star products.  Then $\star$ and
    $\star'$ are $\lie{g}$-equivariantly (strongly) Morita equivalent
    if and only if there is a $\lie{g}$-equivariant symplectomorphism
    $\Psi$ such that $\Psi^*c^{\lie{g}}(\star') - c^{\lie{g}}(\star)$
    is in the image of the first map in 
    \begin{equation}
        \label{eq:EquivariantToInvariantdRhcombined}
        \mathrm{H}^2_{\lie{g}}(M, \mathbb{C})
        \longrightarrow
        \HdR^2(M, \mathbb{C})^{\lie{g}}
        \longrightarrow
        \HdR^2(M,\mathbb{C}),
    \end{equation} 
    and maps to a $2\pi\I$-integral deRham cohomology class under the
    second map.
\end{theorem}
The proof of this theorem is based on a careful analysis of the
semiclassical limit of $\lie{g}$-equivariant equivalence bimodules on
one hand and the invariant Fedosov construction on the other hand. For
star products with quantum momentum maps the proof simplifies in so
far as the problem of lifting the Lie algebra action from the algebras
to the bimodule is trivial thanks to the surjectivity of
\eqref{eq:PicHtoPicForgetAllStuff}. This brings together both parts of
the paper.

Let us finally mention several open questions and further possible
developments: On one hand, the case of Poisson manifolds is clearly
challenging. Here one can expect that things go more or less the same
way. However, it will require more advanced technology: first a
classification of invariant star products has to be achieved. Then the
characterization of Morita equivalent star products
\cite{jurco.schupp.wess:2002a, bursztyn.dolgushev.waldmann:2009a:pre}
has to be cast into the equivariant setting. For both steps the
equivariant formality theorems based on the existence of a
$\lie{g}$-invariant connection \cite{dolgushev:2005a} should be
sufficient. On the other hand, one can pass to higher rank bundles and
hence to endomorphism algebras $\Schnitte(\End(E))$ and their
deformations. This would eventually lead to a notion of equivariant
cohomology with values in endomorphism bundles. From our general
surjectivity statements in \eqref{eq:PicHtoPicForgetAllStuff} one can
expect to get also here some nice applications for the geometric
lifting problems.

The paper is organised as follows: in Section~\ref{sec:Preliminary} we
recall the basic definitions from equivariant Morita theory and
explain the kernels of the groupoid morphisms
\eqref{eq:PicHtoPicForgetAllStuff}. Section~\ref{subsec:FedosovStarProducts}
contains back-ground information on deformation quantization, in
particular on the $\lie{g}$-invariant star products and their
classification where we also reformulate the well-known criteria on
the existence of quantum momentum maps in terms of equivariant
cohomology. In Section~\ref{sec:MomemtumMapsSplitting} we discuss the
general algebraic framework of equivariant Morita theory for algebras
with momentum maps and we proof the first main result of this paper,
the surjectivity of \eqref{eq:PicHtoPicForgetAllStuff}, in
Theorem~\ref{theorem:splitting the image starpich to starpic}.
Finally, Section~\ref{sec:EMEStarProducts} contains the discussion of
the semiclassical limit of equivariant equivalence bimodules of star
products as well as the proof of Theorem~\ref{theorem:TheRealStuff}
both in the general case and with quantum momentum maps.

%
%

\section{Preliminaries on equivariant Morita theory}
\label{sec:Preliminary}

We briefly recall the most important notions of equivariant Morita
equivalence for unital rings and for $^\ast$-algebras over ordered
rings essentially following \cite{jansen.waldmann:2006a}. As a first
remark we note that the terms ``covariant'' and ``equivariant'' are
used synonymously in this context.

%
%

\subsection{Picard groupoids}
\label{subsec:PicardGroupoids}

In the following, we will meet three types of Morita theory: the
original one based on rings alone and two more specific notions taking
into account additional structures on the rings in question. In
particular, we are interested in $^*$-algebras over rings $\ring{C}$
where $\ring{C} = \ring{R}(\I)$ is a ring extension of an ordered ring
$\ring{R}$ by a square root $\I$ of $-1$, see
\cite{bursztyn.waldmann:2005b} for a detailed discussion. The main
example for us, beside $\mathbb{R}$ and $\mathbb{C}$, will be the ring
$\mathbb{R}[[\lambda]]$ of real formal power series with its natural
ordering and the corresponding $\mathbb{C}[[\lambda]]$.

In the case of $^*$-algebras the theory works well if we add the
requirement that all algebras are idempotent and
non-degenerate. However, to keep things simple, we assume that all
rings and algebras are \emph{unital} throughout this work.

We base our approach to Morita theory on ``invertible bimodules''. A
$(\mathcal{B}, \mathcal{A})$-bimodule $\BEA$ is called
\emph{invertible} if one can find another bimodule $\AtildeEB$ such
that $\BEA \tensor[\mathcal{A}] \AtildeEB \cong \BBB$ and $\AtildeEB
\tensor[\mathcal{B}] \BEA \cong \AAA$ as bimodules.  In fact, one
considers isomorphism classes of $(\mathcal{B},
\mathcal{A})$-bimodules as generalized arrows from $\mathcal{A}$ to
$\mathcal{B}$ with composition given by the tensor product and units
given by the isomorphism classes of the canonical bimodules $\AAA$. As
we are dealing with unital rings and algebras throughout, we require
$\Unit_{\mathcal{B}} \cdot x = x = x \cdot \Unit_{\mathcal{A}}$ for
all the elements $x$ in every $(\mathcal{B}, \mathcal{A})$-bimodule.
This way, one obtains the category $\BiMod$ and the invertible
bimodules correspond precisely to the invertible morphisms in this
category which constitute the Picard groupoid $\Pic$. Strictly
speaking, one should fix a suitable Grothendieck universe to work
within. We will assume that this issue has been settled.  It is a
classical theorem of Morita that invertible bimodules are
characterized to be finitely generated full projective right
$\mathcal{A}$-modules with $\mathcal{B} \cong
\End_{\mathcal{A}}(\EA)$. Invertible bimodules are also called
\emph{equivalence bimodules} and rings allowing for an equivalence
bimodule are called \emph{Morita equivalent}.

In the case of $^*$-algebras we require more structure on the
bimodules: first we always assume to have a compatible
$\ring{C}$-module structure on all occurring bimodules.  For $\BEA$ we
require an $\mathcal{A}$-valued inner product, i.e. a map $\SPA{\cdot,
  \cdot}: \BEA \times \BEA \longrightarrow \mathcal{A}$ which is
$\ring{C}$-linear and $\mathcal{A}$-linear to the right in the second
argument, obeys $(\SPA{x, y})^* = \SPA{y, x}$ for all $x, y \in \BEA$,
and which is non-degenerate in the usual sense that $\SPA{x, y} = 0$
for all $y$ implies $x = 0$. Finally, the left multiplications with
elements $b \in \mathcal{B}$ should satisfy $\SPA{b \cdot x, y} =
\SPA{x, b^* \cdot y}$. In particular, the left multiplications are
\emph{adjointable}. Such bimodules are now called \emph{inner product
  bimodules}.  Morphisms between such bimodules are $(\mathcal{B},
\mathcal{A})$-bimodule morphisms which are adjointable. In this
context, isomorphisms always refer to \emph{isometric} (and hence
adjointable) isomorphisms with respect to the inner product.

The canonical bimodule $\AAA$ carries the canonical inner product
$\SP{a, b} = a^*b$. On tensor products we use Rieffel's formula for
an inner product: Let $\CFB$ and $\BEA$ be inner product bimodules then
\begin{equation}
    \label{eq:RieffelInnerProduct}
    \SPFEA{y \tensor[\mathcal{B}] x, y' \tensor[\mathcal{B}] x'}
    =
    \SPEA{x, \SPFB{y, y'} \cdot x'}
\end{equation}
for $y, y' \in \CFB$ and $x, x' \in \BEA$ extends to a well-defined
inner product on $\CFB \tensor[\mathcal{B}] \BEA$ obeying all
requirements except that it may be degenerate. Fortunately, the
degeneracy space is a $(\mathcal{C}, \mathcal{A})$-subbimodule. Thus
we quotient by it and get a new inner product $(\mathcal{C},
\mathcal{A})$-bimodule denoted by $\CFB \tensM[\mathcal{B}] \BEA$.
With the (unital) $^*$-algebras as objects and isomorphism classes of
inner product bimodules as morphisms, $\AAA$ as units, and the
internal tensor product $\tensM[\mathcal{B}]$ as composition, one
obtains the category $\starBiMod$ in analogy to $\BiMod$. The
invertible bimodules in this category are now called
\emph{$^*$-equivalence bimodules} and two $^*$-algebras allowing for
such a $^*$-equivalence bimodule are called \emph{$^*$-Morita
  equivalent}, see \cite{ara:1999a}. The resulting groupoid of
invertible arrows in $\starBiMod$ is called the $^*$-Picard groupoid
$\starPic$.

Up to now we have not yet used the underlying order structure of
$\ring{R}$: Recall that a $\ring{C}$-linear functional $\omega:
\mathcal{A} \longrightarrow \ring{C}$ is called \emph{positive} if
$\omega(a^*a) \ge 0$ for all $a \in \mathcal{A}$ in the sense of
$\ring{R} \subseteq \ring{C}$.  An algebra element $a \in \mathcal{A}$
is called \emph{positive} if for all positive $\omega$ we have
$\omega(a) \ge 0$. This allows for the following definition: an inner
product $\SPA{\cdot, \cdot}$ is called \emph{completely positive} if
for all $n$ and all $x_1, \ldots, x_n \in \BEA$ the matrix $(\SPA{x_i,
  x_j}) \in M_n(\mathcal{A})$ is positive. In this case, $\BEA$ is
also called a \emph{pre Hilbert $(\mathcal{B},
  \mathcal{A})$-bimodule}.  As a matter of fact, the internal tensor
product of two completely positive inner products is again completely
positive. This way, one obtains the subcategory of pre Hilbert
bimodules $\StrBiMod$ of $\starBiMod$. The invertible arrows are now
called \emph{strong equivalence bimodules} and two algebras allowing
for a strong equivalence bimodule are called \emph{strongly Morita
  equivalent}, see \cite{bursztyn.waldmann:2005b}.  The invertible
arrows constitute now the \emph{strong Picard groupoid} denoted by
$\StrPic$. Finally, we note that for $^*$-equivalence bimodules the
inner product \eqref{eq:RieffelInnerProduct} is already
non-degenerate. Hence we do not need the quotient procedure here.

The invertible inner product bimodules in $\starBiMod$ can now be
characterized as follows: they are equivalence bimodules and there is
a $\mathcal{B}$-valued inner product $\BSP{\cdot, \cdot}$ with
$\mathcal{B}$-linearity properties to the left and compatibility with
the right $\mathcal{A}$-multiplications such that both inner products
are \emph{full} in the sense that $\SPA{\mathcal{E}, \mathcal{E}} =
\mathcal{A}$ and $\BSP{\mathcal{E}, \mathcal{E}} = \mathcal{B}$, and
such that they are compatible in the sense that $\BSP{x, y} \cdot z =
x \cdot \SPA{y, z}$ for all $x, y, z \in \BEA$. In the strong case,
$\BSP{\cdot, \cdot}$ is in addition completely positive as well. One
can then show that a candidate for the inverse of $\BEA$ is the
\emph{complex conjugate bimodule} $\AccEB$ with its conjugate
structures, see \cite{bursztyn.waldmann:2005b}.

%
%

\subsection{Equivariant Picard groupoids}
\label{subsec:EquivariantPicardGroupouids}

In a next step we summarize some notions and results concerning
\emph{equivariant} (or \emph{covariant}) Morita theory following
\cite{jansen.waldmann:2006a}.  In the case where the algebras under
investigation possess some group or Lie algebra action, it is natural
to ask for an equivariant notion of equivalence. In order to deal with
groups and Lie algebras simultaneously, we adopt the more general
notion of Hopf algebra actions. In the following, $H$ will denote a
Hopf algebra over $\mathsf{C}$ with comultiplication $\Delta$, counit
$\epsilon$, and antipode $S$. For the comultiplication we use
Sweedler's notation $\Delta(g) = g_\sweedler{1} \tensor
g_\sweedler{2}$.  Moreover, for the $^{\ast}$-algebra case we will
need in addition a $^*$-involution for the Hopf algebra itself such
that the coproduct and the counit are $^*$-homomorphisms, i.e
$(g_\sweedler{1})^* \tensor (g_\sweedler{2})^* = (g^*)_\sweedler{1}
\tensor (g^*)_\sweedler{2}$ and $\epsilon(g^\ast)
=\cc{\epsilon(g)}$. It follows that $S(S(g^\ast)^\ast) = g$ for all $g
\in H$. In particular, $S$ is invertible with $S^{-1} = {}^* S {}^*$.

A \emph{left action} of a Hopf algebra $\acts: H \otimes \alg{A}
\longrightarrow \alg{A}$ on an algebra $\alg{A}$ is then a
$\ring{C}$-linear mapping with $(gh) \acts a = g \acts (h \acts a)$,
$h\acts (ab) = (h_\sweedler{1} \acts a)(h_\sweedler{2}\acts b)$,
$\Unit_H \acts a = a$, and $h\acts \Unit_{\alg{A}} = \epsilon(h)
\Unit_{\alg{A}}$ for all $g, h \in H$ and $a, b \in \alg{A}$. In this
case $\alg{A}$ is also called a \emph{left $H$-module algebra}.  If in
addition $\alg{A}$ and $H$ possess $^{\ast}$-structures we demand
compatibility such that $ (h \acts a)^{\ast} = S(h)^{\ast} \acts
a^{\ast}$ for all $a \in \alg{A}$ and $h \in H$. In this case, the
action will be called a \emph{left $^\ast$-action}.

We now turn to $H$-equivariant representations of left $H$-module
$^{\ast}$-algebras. An $H$-equivariant bimodule $\BEA$ is a bimodule
$\BEA$ with an $H$-action such that $h \acts (b \cdot
x)=(h_\sweedler{1} \acts b)\cdot (h_\sweedler{2} \acts x)$ and $h
\acts (x \cdot a)=( h_\sweedler{1} \acts x) \cdot (h_\sweedler{2}
\acts a)$. The action is called compatible with given $\alg{A}$- or
$\alg{B}$-valued inner products if
\begin{equation}
    \label{eq:2kompatiblilaet H-wirkung}
    h \acts \SPA{x,y}
    =
    \SPA{S(h_\sweedler{1} )^{\ast}\acts x,h_\sweedler{2}\acts y}
    \quad
    \textrm{or}
    \quad
    h \acts \BSP{x, y}
    =
    \BSP{h_\sweedler{1} \acts x, S(h_\sweedler{2})^\ast\acts y},
\end{equation}
respectively.  If $\BEA$ is even a ($^{\ast}$-, strong) equivalence
bimodule with compatible $H$-action, then $\BEA$ is called a
\emph{$H$-equivariant ($^{\ast}$-, strong) equivalence bimodule} and
$\alg{A}$ and $\alg{B}$ are called \emph{$H$-equivariantly
  ($^{\ast}$-, strongly) Morita equivalent}, respectively, see
\cite{jansen.waldmann:2006a}.

There is a canonical action on the internal tensor product of two
equivariant bimodules $\CFB$ and $\BEA$, namely $h \acts (x
\tensor[\alg{B}]{y}) = (h_\sweedler{1} \acts x)
\tensor[\alg{B}](h_\sweedler{2}\acts y)$.  If the bimodules carry
compatible inner products, this action is compatible with the
degeneration space of \eqref{eq:RieffelInnerProduct}. Hence we can
pass to the quotient in this case as well. In addition, given an
$H$-equivariant bimodule $\BEA$ with action $\acts$, there is a
canonical action $\ccacts$ on the complex conjugate $\AccEB$, defined
by $h \ccacts \cc{x}= \cc{S(h)^{*}\acts x}$. All possible
compatibilities of $\acts$ translate to $\ccacts$.

In all three situations (ring-theoretic, $^*$-, and strong) we obtain
a well-defined composition yielding categories $\BiModH$,
$\starBiModH$, and $\StrBiModH$, respectively, such that the
equivalence bimodules correspond to invertible arrows in these
categories. This determines the corresponding equivariant Picard
groupoids, see \cite[Def.~4.1]{jansen.waldmann:2006a}:
\begin{definition}[$H$-equivariant Picard groupoid]
    The unital ($^{\ast}$-) algebras over $\ring{C}$ with compatible
    ($^*$-) action of $H$ as objects and the equivalence classes of
    $H$-equivariant ($^{\ast}$-, or strong) equivalence bimodules as
    morphisms form the $H$-equivariant ($^*$-, or strong) Picard
    groupoids denoted by $\PicH$ ($\starPicH$, or $\StrPicH$).  The
    equivalence classes are constructed with $H$-equivariant
    (isometric) bimodule isomorphisms. The composition of the
    bimodules is the internal tensor product.
\end{definition}

%
%

\subsection{The groups $\Glha$ and $\Uha$}
\label{subsec:GroupsGlhaAndUha}

Forgetting the $H$-action yields natural groupoid morphism
\begin{equation}
    \label{eq:NaturalForgetting}
    \PicH \longrightarrow \Pic,
    \quad
    \starPicH \longrightarrow \starPic,
    \quad
    \textrm{and}
    \quad
    \StrPicH \longrightarrow \StrPic.
\end{equation}
We shall now recall the results from
\cite[Sec.~4.2]{jansen.waldmann:2006a} where the \emph{kernels} of
these groupoid morphisms were studied.  The basic ingredients for the
study of the kernels are the following groups as introduced in
\cite[Def.~A.1]{jansen.waldmann:2006a}. As we will need them later on
we give some details of the construction.
\begin{definition}[The groups $\Glha$ and $\Uha$)]
    \label{definition:GroupsGlHAUHA}
    Let $\alg{A}$ be an $H$-module algebra.  We call $\Glha$ the set
    of all elements $\mathsf{a} \in \Hom_{\mathsf{C}}(H,\alg{A})$ such
    that for all $g, h \in H$ and $b\in \alg{A}$
    \begin{definitionlist}
    \item \label{item:Unital} $\mathsf{a}(\Unit_H) = \Unit_{\alg{A}}$,
    \item \label{item:Product} $\mathsf{a}(gh) = \mathsf{a}
        (g_\sweedler{1}) (g_\sweedler{2}\acts \mathsf{a}(h))$,
    \item \label{item:Commutator} $(h_\sweedler{1}\acts
        b)\mathsf{a}(h_\sweedler{2}) =
        \mathsf{a}(h_\sweedler{1})(h_\sweedler{2}\acts b)$.
    \end{definitionlist}
    If in addition $H$ is a Hopf $^{\ast}$-algebra, $\mathcal{A}$ a
    $^*$-algebra and $\acts$ is a $^*$-action, we call $\Uha$ the set
    of all those elements of $\Glha$ such that for all $h \in H$
    \begin{definitionlist}
        \addtocounter{enumi}{3}
    \item \label{item:Unitary}
        $\mathsf{a}(h_\sweedler{1})(\mathsf{a}(S(h_\sweedler{2})^*))^*
        = \epsilon(h)\Unit_{\alg{A}}$.
    \end{definitionlist}
\end{definition}

As has been shown in \cite[Prop.~A.2]{jansen.waldmann:2006a}, these
sets become groups with respect to the usual convolution product
\begin{equation}
    \label{eq:Konvolutionsprod}
    (\mathsf{a} \ast \mathsf{b})(h)
    =
    \mathsf{a}(h\sweedler{1} )\mathsf{b}(h\sweedler{2})
\end{equation} 
of $\Hom_{\ring{C}}(H, \mathcal{A})$. The unit $\mathsf{e}$ is given
by $\mathsf{e}(g) = \epsilon(g) \Unit_{\mathcal{A}}$. The inverse of
$\mathsf{a}$ is given by $\mathsf{a}^{-1}(g) = g_\sweedler{1} \acts
\mathsf{a}(S(g_\sweedler{2}))$. Moreover, in the $^*$-algebra case,
$\Hom_{\ring{C}}(H, \mathcal{A})$ is a unital $^*$-algebra with
respect to the convolution product and the $^*$-involution defined by
\begin{equation}
    \label{eq:ConvolutionInvolution}
    \mathsf{a}^* (g) = \mathsf{a}(S(g)^*)^*
\end{equation}
for $g \in H$ and $\mathsf{a} \in \Hom_{\ring{C}}(H,
\mathcal{A})$. Then the last condition \refitem{item:Unitary} can be
interpreted as unitarity in the sense that it means $\mathsf{a}^{-1} =
\mathsf{a}^*$ for $\mathsf{a} \in \Uha$.

For later use we mention the following almost trivial situation where
the $^*$-algebra $\mathcal{A}$ is commutative and the action is the
trivial action.
\begin{lemma}
    \label{lemma:TrivialActionOnCommutativeAlgebra}%
    Let $\mathcal{A}$ be a commutative unital $^*$-algebra (algebra)
    and let the $^*$-action (action) of $H$ be trivial, i.e. $g \acts
    a = \epsilon(g)a$ for all $a \in \mathcal{A}$ and $g \in H$. Then
    $\Uha$ consists precisely of the unital $^*$-homomorphisms
    $\mathsf{a}: H \longrightarrow \mathcal{A}$. In the ring-theoretic
    setting $\Glha$ consists of the unital homomorphisms,
    respectively.
\end{lemma}

The relevance of these groups is now the following.  Given an
$H$-action $\acts$ on a Morita equivalence bimodule $\BEA$, we can
construct another action $\acts^{\mathsf{b}}$ on $\BEA$ for all
$\mathsf{b} \in \Glhb$ by
\begin{equation}
    \label{eq:Construction of action of Glha}
    h \acts^{\mathsf{b}} x
    = \mathsf{b}(h_\sweedler{1}) \cdot (h_\sweedler{2} \acts x).
\end{equation}
Indeed, the defining properties of $\mathsf{b}$ guarantee that this is
again an action compatible with the bimodule structure. In the case of
$^*$-algebras and inner product bimodules with compatible actions,
$\mathsf{b} \in \Uhb$ yields an action again compatible with the inner
products. On the other hand, if two $H$-actions $\acts$ and
$\actsp$ are given,
\begin{equation}
    \label{eq:two actions lead to glha}
    u_h(x) 
    = h_\sweedler{1}\acts (S(h_\sweedler{2})\actsp x)
\end{equation}
for $h \in H$ and $ x \in \BEA$ defines a right $\alg{A}$-linear
endomorphism $u_h$ of $\BEA$. Because we work in the unital case,
$u_h$ is the left multiplication by a unique element of $\alg{B}$,
which we denote by $\mathsf{b}(h)$. It turns out that $\mathsf{b} \in
\Glhb$ or, for $^*$-equivalence bimodules, $\mathsf{b} \in \Uhb$ if
the actions were compatible with the inner products. It is now easy to
see that this gives indeed a bijection and hence the elements of
$\Glhb$ (or $\Uhb$) parametrize the possible compatible $H$-actions on
a given equivalence bimodule (or $^*$-equivalence bimodule,
respectively), \emph{provided} there is at least one compatible action
at all.  In fact, the above construction gives a free and transitive
group action of $\Glhb$ (or $\Uhb$) on the set of possible compatible
$H$-actions.

Though $\Glhb$ (or $\Uhb$, respectively) acts freely and transitively
on the set of all compatible $H$-actions, some of them might lead to
isomorphic bimodules and hence to the same element in the Picard
groupoid. They can be described as follows: We denote by $\Glzb$ the
group of all central and invertible elements of $\alg{B}$ and by
$\Uzb$ the group of all central and unitary elements, respectively.
The subgroups of $H$-invariant elements are denoted by $\Glzb^H$ and
$\Uzb^H$, respectively.  There is a group morphism
\begin{equation}
    \label{eq:von GlZ(A) nach GL(H,A)}
    \Glzb \ni c
    \longmapsto
    \left(
        \hat{c}:h \longmapsto c(h \acts c^{-1})
    \right)
    \in \Glhb.
\end{equation}
Then $\acts$ and $\acts^{\mathsf{b}}$ define isomorphic
$H$-equivariant equivalence bimodules iff $\mathsf{b}$ is in the image
of \eqref{eq:von GlZ(A) nach GL(H,A)}. The analogous result holds for
the case of $^*$-algebras and the groups $\Uzb$ and $\Uhb$ instead of
$\Glzb$ and $\Glhb$. It follows that the image of \eqref{eq:von GlZ(A)
  nach GL(H,A)} is a normal subgroup such that the quotient groups
$\Glohb =\Glhb / \widehat{\Glzb}$ and $\Uohb =\Uhb / \widehat{\Uzb}$
parametrize the \emph{inequivalent} actions on an equivalence
bimodule. In particular, we obtain short exact sequences
\begin{equation}
    \label{eq:kernel picH-pic}
    \xy
    \morphism(0,0)<700,0>[1`\Glohb;]
    \morphism(700,0)<900,0>[\Glohb`\PicH(\alg{B});]
    \morphism(1600,0)<1000,0>[\PicH(\alg{B})`\Pic(\alg{B});]
    \endxy
\end{equation}
in the ring-theoretic situation and
\begin{equation}
    \label{eq:kernel starpicH-starpic}
    \xy
    \morphism(0,0)<700,0>[1`\Uohb;]
    \morphism(700,0)<900,0>[\Uohb`\starPicH(\alg{B});]
    \morphism(1600,0)<1000,0>[\starPicH(\alg{B})`\starPic(\alg{B});]
    \endxy
\end{equation}  
in the $^*$-algebra framework, thereby encoding the kernels of the
groupoid morphisms $\PicH \longrightarrow \Pic$ and $\starPicH
\longrightarrow \starPic$. Since the complete positivity of inner
products was not playing any role in the above argument, the analogous
statement also holds for the strong version.

Similarly, there is an action of $\mathsf{a} \in \Glha$ (or
$\mathsf{a} \in \Uha$, respectively) from the right on $(\BEA, \acts)$
by twisting the $H$-action $\acts$ to $g \acts_{\mathsf{a}} x =
(g_\sweedler{1} \acts x) \cdot \mathsf{a}(g_\sweedler{2})$ for $g \in
H$ and $x \in \BEA$. This can be used to see that the groups $\Glhb$
and $\Glha$ are actually isomorphic for $H$-equivariantly Morita
equivalent algebras.  Analogously $\Uhb$ and $\Uha$ are isomorphic for
$H$-equivariantly $^*$-Morita equivalent $^*$-algebra. For later use
we recall the relation between the two groups for the $^*$-algebra
case:
\begin{lemma}
    \label{lemma:Isomorphism uha uhb}
    Let $\BEA \in \starPicH(\alg{B}, \alg{A})$ and $\mathsf{b} \in
    \Uhb$ and $\mathsf{a}\in \Uha$. We have
    $(\acts_{\mathsf{a}})^{\mathsf{b}}=(\acts^{\mathsf{b}})_{\mathsf{a}}$
    and there exists a unique $h_{\mathcal{E}}(\mathsf{a})\in \Uhb$,
    such that
    $\acts_{\mathsf{a}}=\acts^{h_{\mathcal{E}}(\mathsf{a})}$.  The map
    \begin{equation}
        \label{eq:uha nach uhb}
        h_{\mathcal{E}}: \Uha \ni \mathsf{a} \longmapsto
        h_{\mathcal{E}}(\mathsf{a}) \in \Uhb
    \end{equation}
    is a group isomorphism and determines a group isomorphism $\Uoha
    \longrightarrow \Uohb$ that will also be denoted by
    $h_{\mathcal{E}}$.
\end{lemma}
In fact, $h_{\mathcal{E}}$ turns out to come from a groupoid action of
$\starPicH$ on the collection of all the groups $\mathsf{U}(H, \cdot)$
in a very precise sense, see \cite[Lem.~5.11,
Thm~5.14]{jansen.waldmann:2006a} for details and the proof of this
lemma.

%
%

\section{Fedosov star products and quantum momentum maps}
\label{subsec:FedosovStarProducts}

While the above considerations apply to general algebras and
$^*$-algebras, respectively, we focus now on the star product algebras
arising from deformation quantization \cite{bayen.et.al:1978a}. To
keep things simple, we consider star products on a \emph{symplectic}
manifold which we assume to be connected. We recall the basic notions.

%
%

\subsection{Deformation quantization}
\label{subsec:DeformationQuantization}

A \emph{star product} $\star$ is a $\mathbb{C}[[\lambda]]$-bilinear
associative multiplication for the formal power series with
coefficients in the complex-valued functions $C^\infty(M)[[\lambda]]$
written as
\begin{equation}
    \label{eq:StarProduct}
    f \star g = \sum_{r=0}^\infty \lambda^r C_r(f, g)
\end{equation}
for $f, g \in C^\infty(M)[[\lambda]]$ such that $C_0(f, g) = fg$ is
the pointwise (commutative) product and $C_1(f, g) - C_1(g, f) = \I
\{f, g\}$ yields the Poisson bracket. Moreover, we require $\star$ to
be \emph{differential}, i.e. all the maps $C_r$ should be
bidifferential operators on $M$. Finally, one requires $1 \star f = f
= f \star 1$.  A star product is called \emph{Hermitian} if in
addition $\cc{f \star g} = \cc{g} \star \cc{f}$. Here the complex
conjugation of $\lambda$ is defined to be $\cc{\lambda} =
\lambda$. Thus a Hermitian star product turns $C^\infty(M)[[\lambda]]$
into a unital $^*$-algebra over $\mathbb{C}[[\lambda]]$ and we arrive
at the general framework discussed before. Finally, recall that two
star products $\star$ and $\star'$ are called \emph{equivalent} if
there is a formal series $T = \id + \sum_{r=1}^\infty \lambda^r T_r$
of differential operators such that $f \star' g = T^{-1} (Tf \star
Tg)$, i.e. $T$ is an algebra isomorphism deforming the identity. On a
symplectic manifold, there is a characteristic class which assigns to
every $\star$ a formal series $c(\star) \in \frac{[\omega]}{\I
  \lambda} + \HdR^2 (M, \mathbb{C})[[\lambda]]$ in the second deRham
cohomology such that two star products are equivalent iff their
characteristic classes coincide, see e.g. \cite{gutt.rawnsley:1999a}
for a nice introduction. Every such formal series arises as the
characteristic class of a star product.  Note that the characteristic
class is natural with respect to symplectomorphisms,
i.e. $\Psi^*c(\star) = c(\Psi^*(\star))$, when we move around the star
product by means of $\Psi$. A pedagogical introduction to deformation
quantization and further references can be found in
\cite{waldmann:2007a}.

In the following we shall make extensive use of the Fedosov
construction \cite{fedosov:1994a} of star products on a symplectic
manifold. Fedosov constructed out of a symplectic torsion-free
connection $\nabla$ and a formal series of closed two forms $\Omega
\in \lambda \Omega^2(M)[[\lambda]]$ a differential star product, the
\emph{Fedosov star product} $\star_{\nabla, \Omega}$. It is Hermitian
iff $\cc{\Omega} = \Omega$ is real. It turns out that
$c(\star_{\nabla, \Omega}) = \frac{1}{\I\lambda} [\omega + \Omega]$,
see e.g. \cite[Sect.~6.4]{waldmann:2007a} for further details.  We
shall also need the following simple generalization of Fedosov's
construction \cite{waldmann:2002b}: let $E \longrightarrow M$ be a
complex vector bundle over a symplectic manifold $M$ and choose a
covariant derivative $\nabla^E$ for it. Then one can construct a
formal series of bidifferential operators $R_r: \Schnitte(E) \times
C^\infty(M) \longrightarrow \Schnitte(E)$ such that
\begin{equation}
    \label{eq:RightModuleStructure}
    s \bullet f = s \cdot f + \sum_{r=1}^\infty \lambda^r R_r(s, f)
\end{equation}
for $s \in \Schnitte(E)[[\lambda]]$ and $f \in C^\infty(M)[[\lambda]]$
defines a deformation of the classical right module structure into a
right module structure with respect to $\star_{\nabla,
  \Omega}$. Moreover, there exists a formal series of bidifferential
operators $C_r': \Schnitte(\End(E)) \times \Schnitte(\End(E))
\longrightarrow \Schnitte(\End(E))$ such that
\begin{equation}
    \label{eq:AstarprimeB}
    A \star' B = AB + \sum_{r=1}^\infty \lambda^r C_r'(A, B)
\end{equation}
for $A, B \in \Schnitte(\End(E))[[\lambda]]$ defines a formal
associative deformation. Finally, there is also a bidifferential
deformed left module structure
\begin{equation}
    \label{eq:Abulletprimes}
    A \bullet' s = As + \sum_{r=1}^\infty \lambda^r R'_r (A, s)
\end{equation}
for $A \in \Schnitte(\End(E))[[\lambda]]$ and $s \in
\Schnitte(E)[[\lambda]]$ with respect to $\star'$ such that in total
we obtain a deformed \emph{bimodule} structure. In fact, this turns
out to be a ring-theoretic equivalence bimodule and all equivalence
bimodules with $C^\infty(M)[[\lambda]]$ acting from the right are
isomorphic to such a deformed vector bundle. For the classical limit
this is of course well-known. In the case $E = L$ is a line bundle,
$\star'$ yields a star product, too, and the deformed line bundle is
an equivalence bimodule for $\star$ and $\star'$.

In addition, one can specify a pseudo Hermitian fiber metric $h$ on
$E$, viewed as an inner product $h: \Schnitte(E) \times \Schnitte(E)
\longrightarrow C^\infty(M)$. This endows $\Schnitte(\End(E))$ with a
$^*$-involution. Next, if we assume that $\nabla^E$ is compatible with
$h$, the Fedosov construction automatically yields a $^*$-equivalence
bimodule by suitably deforming $h$ as well. If in addition $h$ was a
Hermitian fiber metric then the deformed $^*$-equivalence bimodule is
even a strong one. Again, all are of this form up to isomorphisms.

%
%

\subsection{Quantum momentum maps}
\label{subsec:QuantumMomentumMaps}

We will now start to incorporate symmetries into this situation.  In
the following we will always assume that there is a Lie algebra
$\lie{g}$ acting on $(M, \omega)$ by symplectic vector fields denoted
by $X_\xi$ where $\xi \in \lie{g}$. We follow the convention that
$\rho(\xi) = - \Lie_{X_\xi}$ defines a left representation of
$\lie{g}$ on $C^\infty(M)$, where $\Lie$ denotes the usual Lie
derivative.  Thus $\rho(\xi)$ is a Poisson derivation. The typical
situation we have in mind is that the $X_\xi$ are the fundamental
vector fields of a symplectic Lie group action. However, in the
following we never need the integrated version but exclusively use the
Lie algebraic point of view. This opens the framework also to more
general Lie algebra actions where the fundamental vector fields may be
non-complete. However, throughout this section, we do assume that on
$M$ there is a $\lie{g}$-invariant connection. In this case, there is
also a $\lie{g}$-invariant symplectic torsion-free connection
$\nabla$.

To cast this situation into the general Hopf $^*$-algebra framework we
consider the following two Hopf algebras: first, we take the formal
power series in the complexified tensor algebra over $\lie{g}$, viewed
as a $\mathbb{C}[[\lambda]]$-module, and divide by the two-sided ideal
generated by the relations $\xi \otimes \eta - \eta \otimes \xi - \I
\lambda [\xi, \eta]$ for all $\xi, \eta \in \lie{g}$. The resulting
``rescaled'' universal enveloping algebra will be denoted by
$\mathcal{U}_\lambda (\lie{g})$.  The $\I$ in front of the bracket
allows to consider the elements $\xi \in \lie{g} \subseteq
\mathcal{U}_\lambda(\lie{g})$ to be Hermitian, the reason for the
additional $\lambda$ will become clear in
\eqref{eq:JisHomomorphism}. The requirements
\begin{equation}
    \label{eq:structures of rescaled universal enveloping}
    \xi^* =  \xi,
    \quad
    \Delta(\xi)= \xi \otimes \Unit + \Unit \otimes \xi,
    \quad
    \epsilon(\xi) = 0,
    \quad
    \textrm{and}
    \quad
    S(\xi) = -\xi
\end{equation}
for all $\xi \in \lie{g}$ lead to a unique Hopf $^*$-algebra structure
for $\mathcal{U}_\lambda(\lie{g})$ over
$\mathbb{C}[[\lambda]]$. Second, the complexified universal enveloping
algebra $\mathcal{U}_{\mathbb{C}}(\lie{g})$ is a Hopf $^{*}$-algebra,
too, with $^{*}$-involution determined by $\xi^{*} = - \xi$. To be in
the framework of Hopf $^*$-algebras over $\mathbb{C}[[\lambda]]$ we
extend all structure maps to
$\mathcal{U}_{\mathbb{C}}(\lie{g})[[\lambda]]$.

We have to strengthen now the assumptions about the action in the
following way: Recall that, on the classical level, the Lie algebra
action is called \emph{Hamiltonian} if there is a linear mapping $J_0:
\lie{g} \longrightarrow C^\infty(M)$ with $\ins_{X_{\xi}}\omega = \D
J_0(\xi)$. Equivalently, this means that the Poisson derivations
$\rho(\xi) = \{J_0(\xi), \cdot\}$ are \emph{inner}.  If in addition
$J_0$ is equivariant with respect to the adjoint action of $\lie{g}$,
which is equivalent to
\begin{equation}
    \label{eq:classical momentum map}
    J_0([\xi,\eta]) = \{J_0(\xi), J_0(\eta)\}
    \quad
    \textrm{for all}
    \quad
    \xi, \eta \in \lie{g},
\end{equation}
then $J_0$ is called \emph{classical momentum map}.

The appropriate notion of symmetry in deformation quantization is now
the following, see e.g. \cite{bertelson.bieliavsky.gutt:1998a} and
references therein. A star product $\star$ is called
\emph{$\lie{g}$-invariant} if the action $\rho $ of $\lie{g}$ is given
by derivations of $\star$.  The Fedosov construction depends
functorially on the data $\nabla$, $\Omega$. This allows to show that
$\star_{\nabla, \Omega}$ is $\lie{g}$-invariant iff the entrance data
$\nabla$ and $\Omega$ are $\lie{g}$-invariant. Moreover, one can show
that every $\lie{g}$-invariant star product $\star$ is
$\lie{g}$-equivariantly equivalent to a $\lie{g}$-invariant Fedosov
star product $\star_{\nabla, \Omega}$. Finally, one knows that
$\star_{\nabla, \Omega}$ is $\lie{g}$-equivariantly equivalent to
$\star_{\nabla', \Omega'}$ iff $\Omega$ and $\Omega'$ are cohomologous
in the $\lie{g}$-invariant deRham cohomology, see
\cite{bertelson.bieliavsky.gutt:1998a} for detailed proofs for the
case of a Lie group action: the Lie algebraic case can be done by
completely analogous means. This allows to define the
\emph{$\lie{g}$-invariant characteristic class}
\begin{equation}
    \label{eq:InvariantCharacteristicClass}
    c^{\lie{g}} (\star)
    = \frac{[\omega + \Omega]}{\I\lambda}
    \in
    \frac{[\omega]}{\I\lambda}
    + \HdR^2(M, \mathbb{C})^{\lie{g}}[[\lambda]],
\end{equation}
where $\Omega$ is a formal series of $\lie{g}$-invariant closed
two-forms such that $\star$ is $\lie{g}$-equivariantly equivalent to
$\star_{\nabla, \Omega}$. Note that under the canonical map
\begin{equation}
    \label{eq:CanonicalHdRgToHdR}
    \HdR^2(M, \mathbb{C})^{\lie{g}}
    \longrightarrow \HdR^2(M, \mathbb{C}),
\end{equation}
extended to \eqref{eq:InvariantCharacteristicClass}, the
$\lie{g}$-invariant characteristic class $c^{\lie{g}}(\star)$ is
mapped to the characteristic class $c(\star)$. In particular, an
arbitrary star product $\star$ is equivalent to a $\lie{g}$-invariant
star product iff $c(\star)$ is in the image of
\eqref{eq:CanonicalHdRgToHdR}. Note also that the construction of
$c^{\lie{g}}(\cdot)$ requires the existence of a $\lie{g}$-invariant
connection but does not depend on the particular choice of it.

In the following, we assume that a classical momentum map $J_0$ exists
and is given.  Then the analog of a classical momentum map is a
quantum momentum map, see \cite{mueller-bahns.neumaier:2004a} and
references therein:
\begin{definition}[Quantum momentum map]
    \label{definition:quantum momentum map}
    Let $\star$ be a $\lie{g}$-invariant star product. A linear map $J
    = J_0 + J_+: \lie{g} \longrightarrow \Fkt[[\lambda]]$ with $J_0:
    \lie{g} \longrightarrow \Fkt$ and $J_+: \lie{g} \longrightarrow
    \lambda \Fkt[[\lambda]]$ is called quantum Hamiltonian for the
    action $\rho$ if
    \begin{equation}
        \label{eq:def quantum Hamiltonian}
        \rho(\xi)
        =
        \frac{1}{\I \lambda}\ad_{\star}(J(\xi))
        \quad
        \textrm{for all}
        \quad
        \xi \in \lie{g}.
    \end{equation}
    Moreover, $J$ is called a quantum momentum map if in addition for
    all $\xi, \eta \in \lie{g}$
    \begin{equation}
        \label{eq:def quantum momentum map}
        \I \lambda J([\xi, \eta]) = [J(\xi), J(\eta)]_{\star}.
    \end{equation}
    In the case where $\star$ is Hermitian we require in addition
    $J(\xi)^{\ast} = J(\xi)$.
\end{definition}

Note that the right hand side of \eqref{eq:def quantum momentum map}
is necessarily of order $\lambda$ making the condition
meaningful. Thanks to this and thanks to the definition of the
rescaled universal enveloping algebra, a quantum momentum map uniquely
extends to a unital $\mathbb{C}[[\lambda]]$-linear algebra homomorphism
\begin{equation}
    \label{eq:JisHomomorphism}
    J: \mathcal{U}_\lambda (\lie{g})
    \longrightarrow
    C^\infty(M)[[\lambda]],
\end{equation}
which is even a $^*$-homomorphism as soon as $\star$ is
Hermitian. Moreover, the definition
\begin{equation}
    \label{eq:xiactsf}
    \xi \acts f = \I \lambda \rho(\xi) f
\end{equation}
for $\xi \in \lie{g}$ and $f \in C^\infty(M)[[\lambda]]$ extends to a
Hopf algebra action (even a $^*$-action in the Hermitian case) of
$\mathcal{U}_\lambda(\lie{g})$ on $C^\infty(M)[[\lambda]]$ as soon as
$\star$ is $\lie{g}$-invariant.

Since we assume to have a $\lie{g}$-invariant connection, one even has
necessary and sufficient conditions for the existence of quantum
momentum maps for $\lie{g}$-invariant Fedosov star products.  In fact,
there exists a quantum momentum map for the $\lie{g}$-invariant
Fedosov star product $\star_{\nabla, \Omega}$ if and only if there
exists a cochain $J \in C^1(\lie{g}, C^\infty(M))[[\lambda]]$ such
that
\begin{equation}
    \label{eq:conditions for quantum momentum map}
    \D J(\xi)
    =
    \ins_{X_{\xi}}(\omega +\Omega)
    \quad
    \textrm{and}
    \quad
    (\omega + \Omega)(X_{\xi}, X_{\eta})
    =
    J([\xi, \eta])
    \quad
    \textrm{for all}
    \quad
    \xi, \eta \in \lie{g}.
\end{equation}
Moreover, in this case the quantum momentum map $J$ is determined up
to Chevalley cocycles in $Z^1(\lie{g}, \mathbb{C})[[\lambda]]$. In
case of a Hermitian Fedosov star product, i.e. $\cc{\Omega} = \Omega$,
the real quantum momentum map is determined up to Chevalley cocycles
in $Z^1(\lie{g}, \mathbb{R})[[\lambda]]$, see
\cite{mueller-bahns.neumaier:2004a} for detailed proofs of these
statements.

For any other $\lie{g}$-invariant star product $\star$ we can pass to the
corresponding Fedosov star product $\star_{\nabla, \Omega}$ by means
of a $\lie{g}$-equivariant equivalence transformation thanks to the
above classification, see
\cite[Prop.~4.1]{bertelson.bieliavsky.gutt:1998a}. So if
\begin{equation}
    \label{eq:Equivalence transformation to Fedosov star product}
    f \star_{\nabla, \Omega} g
    =
    T \left(T^{-1}(f) \star T^{-1}(g)\right)
\end{equation}
denotes such a $\lie{g}$-equivariant equivalence transformation $T$,
it is easy to see that $J$ is a quantum momentum map for $\star$ if
and only if $TJ$ is a quantum momentum map for $\star_{\nabla,
  \Omega}$.  In fact, one can transport quantum momentum maps with
equivariant equivalence transformations to any other equivalent star
product as well. Thus the existence of a quantum momentum map is a
property of the $\lie{g}$-invariant characteristic class
$c^{\lie{g}}(\star)$ within the $\lie{g}$-invariant star products and
not just of the star product itself.

%
%

\subsection{Equivariant cohomology}
\label{subsec:EquivariantCohomology}

We will now rephrase the above well-known conditions for the existence
of quantum momentum maps in terms of equivariant cohomology. To this
end, recall the Cartan model for equivariant cohomology, see
e.g. \cite{guillemin.sternberg:1999a} for details.

The \emph{$\lie{g}$-equivariant differential forms} of degree $k$ are
elements of the space
\begin{equation}
    \label{eq:equivariant forms}
    \Omega^k_{\lie{g}}(M, \mathbb{C})
    = \bigoplus_{2i+j=k}
    \left(
        \Sym^i(\lie{g}^{\ast}) \otimes \Omega^j(M, \mathbb{C})
    \right)^{\lie{g}}.
\end{equation}
Here $\Sym^i(\lie{g}^{\ast})$ denotes the (complexified) symmetric
tensors over $\lie{g}^*$ endowed with the inherited coadjoint action
of $\lie{g}$. We usually view them as homogeneous polynomial maps from
$\lie{g}$ to $\mathbb{C}$ of degree $i$. Thus we can view $\alpha \in
\Omega^k_G(M, \mathbb{C})$ as a polynomial subject to the equivariance
condition
\begin{equation}
    \label{eq:explicit equivariance definition}
    \alpha([\xi, \eta]) = \Lie_{X_\xi}\alpha(\eta)
    \quad
    \textrm{for all}
    \quad    
    \xi, \eta \in \lie{g}.
\end{equation}
The differential $\D_{\lie{g}}: \Omega_{\lie{g}}^k(M, \mathbb{C})
\longrightarrow \Omega_{\lie{g}}^{k+1}(M, \mathbb{C})$ is defined to
be
\begin{equation}
    \label{eq:equivariant differential}
    (\D_{\lie{g}} \alpha)(\xi)
    =
    \D(\alpha(\xi)) + \ins(X_{\xi})\alpha(\xi).
\end{equation}
Then the Cartan model of equivariant cohomology is defined to be the
cohomology of this complex, i.e.
\begin{equation}
    \label{eq:equivariant cohomology}
    \mathrm{H}_{\lie{g}}^{\bullet}(M, \mathbb{C})
    = \ker\D_{\lie{g}} / \image \D_{\lie{g}},
\end{equation}
which is graded as $\D_{\lie{g}}$ is compatible with the degree $k$.
In the following we will exclusively work with the Cartan model of
equivariant cohomology and simply call this the \emph{equivariant
  cohomology} of $M$ with respect to the Lie algebra action of
$\lie{g}$. Note that $\D_{\lie{g}}$ commutes with complex conjugation
and hence allows for a real equivariant cohomology as well.
\begin{remark}
    \label{remark:NotUsualStuff}%
    Note that we do not require that the Lie algebra action integrates
    to a Lie group action, contrary to
    \cite{guillemin.sternberg:1999a}. Even if it integrates, the
    corresponding Lie group action is allowed to be quite
    arbitrary. In particular, we do not need it to be proper.  As
    already mentioned, for our purposes it will turn out that the
    existence of the invariant connection is all we need which is of
    course a strictly weaker requirement than a proper action.
\end{remark}

We are mainly interested in the second equivariant cohomology.  In
degree $k = 1$ and $k = 2$, the equivariant forms are
\begin{equation}
    \label{eq:for k=2 equivariant forms}
    \Omega_{\lie{g}}^1 (M, \mathbb{C})
    =
    \Omega^1 (M, \mathbb{C})^{\lie{g}}
    \quad
    \textrm{and}
    \quad
    \Omega_{\lie{g}}^2 (M, \mathbb{C})
    =
    \Omega^2(M, \mathbb{C})^{\lie{g}} 
    \oplus C^1(\lie{g}, C^{\infty}(M))^{\lie{g}}.
\end{equation}
From this we see that projecting to the first component in a closed
equivariant two-form yields a well-defined canonical map
\begin{equation}
    \label{eq:EquivariantToInvariantdRh}
    \mathrm{H}^2_{\lie{g}}(M, \mathbb{C})
    \longrightarrow
    \HdR^2(M, \mathbb{C})^{\lie{g}}.
\end{equation}

Furthermore, a closed two-form $\Omega \in \Omega^2(M)$ can
be extended to a closed equivariant two-form if and only if there is a
$J \in C^1(\lie{g}, C^{\infty}(M))$ such that
\begin{equation}
    \label{eq:concret conditions for second equivariant cohomology}
    \Lie_{X_\xi} \Omega = 0,
    \quad
    \ins_{X_{\xi}} \Omega = \D J(\xi),
    \quad
    \textrm{and}
    \quad
    J([\xi, \eta]) = \Omega(X_{\xi},X_{\eta})
\end{equation}   
for all $\xi, \eta \in \lie{g}$. Using this, we can reformulate the
existence of quantum momentum maps in terms of the equivariant
cohomology.
\begin{proposition}
    \label{proposition:existence of momentum maps and equivariant cohomology}
    Let $\lie{g}$ act on $(M, \omega)$ symplectically with classical
    momentum map $J_0$. Assume furthermore that there exists a
    $\lie{g}$-invariant symplectic connection $\nabla$ for $M$.  Then
    for a $\lie{g}$-invariant star product $\star$ there exists a
    quantum momentum map if and only if $c^{\lie{g}}(\star)$ is in the
    image of \eqref{eq:EquivariantToInvariantdRh}.
\end{proposition}
\begin{proof}
    First assume $\star$ is equal to the Fedosov star product for
    $\Omega$ and $\nabla$. Then the assertion follows from the
    characterization \eqref{eq:conditions for quantum momentum map},
    see \cite{mueller-bahns.neumaier:2004a}. Since every
    $\lie{g}$-invariant star product is $\lie{g}$-invariantly
    equivalent to a Fedosov star product and since
    $c^{\lie{g}}(\star)$ is precisely defined by means of the
    corresponding Fedosov data, the general case follows as well by
    transporting the quantum momentum map via the equivalence.
\end{proof}
The corresponding characterization of the existence of a classical
momentum map is well-known, see e.g. the discussion in
\cite[Chap.~9]{guillemin.sternberg:1999a}.

\begin{example}[The Weyl-Moyal star product]
    \label{example:SPzweinOnRzwein}%
    Let us discuss an illustrative and well-known example for
    invariant star products. For convenience, we formulate this
    example also in the global way using Lie group actions.  We
    consider $\mathbb{R}^{2n}$ with the standard symplectic form
    $\omega$ and the canonical symplectic action of the affine
    symplectic group $\mathrm{Sp}_{2n} \ltimes \mathbb{R}^{2n}$. Note
    that this group action is not proper, the canonical flat
    connection is, however, invariant with respect to this group
    action. The invariant second deRham cohomology
    $\HdR^{2}(\mathbb{R}^{2n}, \mathbb{C})^{\mathrm{Sp}_{2n}\ltimes
      \mathbb{R}^{2n}} \cong \mathbb{C}[\omega]$ is one
    dimensional. In particular, \eqref{eq:CanonicalHdRgToHdR} is
    \emph{not} injective in this case (as it would be for compact
    groups). The map $J_{\mathrm{Sp}}(\xi)(x) = - \frac{1}{2}x^T\omega
    \xi x$ for $\xi$ in the Lie algebra $\mathfrak{sp}_{2n}$ of
    $\mathrm{Sp}_{2n}$ and $x \in \mathbb{R}^{2n}$ defines a classical
    momentum map for the action of $\mathrm{Sp}_{2n}$. In contrast,
    the infinitesimal action of the translations is Hamiltonian with
    respect to the map $J_{\mathbb{R}^{2n}}(a)(x) = - x^T \omega a$
    for $a,x \in \mathbb{R}^{2n}$ but with absence of
    $\ad^{*}$-equivariance.  The Weyl-Moyal star product
    \begin{equation}
        \label{eq:WeylMoyal}
        f \star_0 g =
        \sum_{r=0}^\infty \frac{1}{r!}
        \left(- \frac{\I\lambda}{2}\right)^r
        \omega^{i_1 j_1} \cdots \omega^{i_r j_r}
        \frac{\partial^r f}{\partial x^{i_1} \cdots \partial x^{i_r}}
        \frac{\partial^r g}{\partial x^{j_1} \cdots \partial x^{j_r}},
    \end{equation}
    with $\omega^{ij}$ being the inverse matrix of the coefficient
    matrix of $\omega$, is invariant with respect to the action of
    $\mathrm{Sp}_{2n}\ltimes \mathbb{R}^{2n}$. Moreover,
    $J_{\mathrm{Sp}}(\xi)$ is a quantum momentum map for the
    $\mathrm{Sp}_{2n}$ action, whereas $J_{\mathbb{R}^{2n}}$ yields
    only a quantum Hamiltonian for the translations. By
    $\star_{c\omega}$ for $c \in \lambda \mathbb{C}[[\lambda]]$ we
    denote the ``Weyl-Moyal'' star products for the rescaled
    symplectic form $(1+c)\omega$, i.e. the Fedosov star products to
    the choice $\Omega = c\omega$. It can be shown that these are the
    only $\mathrm{Sp}_{2n}\ltimes \mathbb{R}^{2n}$-invariant star
    products on $\mathbb{R}^{2n}$: this follows from the fact that any
    such invariant star product has to be equivalent to a Fedosov one
    and the fact that there are simply no candidates for
    \emph{invariant} equivalence transformations beside the identity
    due to the lack of invariant differential operators.  Clearly
    $(1+c)J_{\mathrm{Sp}}$ yields a quantum momentum map for
    $\star_{c\omega}$ for the $\mathrm{Sp}_{2n}$-action. In
    particular, the Weyl-Moyal star product $\star_0$ is the
    \emph{unique} $\mathrm{Sp}_{2n}\ltimes \mathbb{R}^{2n}$-invariant
    star product with the classical momentum map for the
    $\mathrm{Sp}_{2n}$-action as quantum momentum map. This feature is
    sometimes referred to as \emph{strong invariance}. The formal
    exponential of the Euler vector field leads to
    $\mathrm{Sp}_{2n}$-invariant equivalences between the different
    star products $\star_{c\omega}$, as expected from the equation
    $[\omega] = [0]$ in $\HdR^2(\mathbb{R}^{2n})^{\mathrm{Sp}_{2n}}$.
\end{example}

%
%

\section{Momentum maps and the splitting}
\label{sec:MomemtumMapsSplitting}

In this section we will generalize the notion of a (quantum) momentum
map to the general algebraic situation of Hopf algebra actions. This
will allow us to prove a general splitting theorem for the groupoid
morphisms \eqref{eq:NaturalForgetting}.

%
%

\subsection{Momentum maps}
\label{subsec:MomentumMaps}

The image of \eqref{eq:NaturalForgetting} is, in general, very hard to
understand and corresponds to a ``lifting'' problem: the $H$-action
has to be lifted from the algebras to the bimodule, see the discussion
in \cite[Ex.~4.8]{jansen.waldmann:2006a}. However, if we restrict to a
subclass of algebras having a momentum map things turn out to be very
nice. The following axiomatization of a (quantum) momentum map is now
straightforward:
\begin{definition}[Momentum map]
     \label{def:momentum map of H}
     Let $\mathcal{A}$ be a unital $^*$-algebra with a $^*$-action of
     a Hopf $^*$-algebra $H$.  Then a momentum map $J$ is a
     $^*$-homomorphism $J: H \longrightarrow \alg{A}$ such that for
     all $h \in H$ and $a \in \mathcal{A}$
     \begin{equation}
         \label{eq:action via momentum map}
         h\acts  a = J(h\sweedler{1} ) a J(S(h\sweedler{2})).
    \end{equation}
    The set of all such momentum maps for the given $^*$-action is
    denoted by $\starMM_H(\mathcal{A})$.
\end{definition}
Conversely, if a $^*$-homomorphism $J: H \longrightarrow \mathcal{A}$
is given, then \eqref{eq:action via momentum map} defines a
$^*$-action of $H$ on $\mathcal{A}$. In the ring-theoretic case we
drop the requirement $J(g^*) = J(g)^*$ and obtain the analogous
definition for a momentum map. Then the set of all momentum maps is
denoted by $\MM_H(\mathcal{A})$, respectively.

Note that the action of $H$ is \emph{inner}, but the choice of the
inner elements is even consistent with the algebra structure of
$H$. This corresponds to the equivariance requirement of classical
momentum maps in symplectic geometry.  Note also that a momentum map
for a non-trivial action of $H$ requires $\mathcal{A}$ to be
sufficiently noncommutative: if the action is inner, the center of
$\mathcal{A}$ consists of $H$-invariant elements, i.e.
\begin{equation}
    \label{eq:CentralElementsInvariant}
    \Zentrum{\mathcal{A}} \subseteq \mathcal{A}^H.
\end{equation}
There is no a priori restriction on the Hopf algebra, as the adjoint
action of any Hopf algebra on itself is just an inner action with the
identity as momentum map.

For a star product algebra $(C^{\infty}(M)[[\lambda]], \star)$, this
notion of a momentum map for $H = \mathcal{U}_{\lambda}(\lie{g})$
coincides with the notion of a quantum momentum map if the action is
of the form \eqref{eq:xiactsf}. Note that it is crucial to use the
\emph{rescaled} universal enveloping algebra since for a star product
algebra \eqref{eq:action via momentum map} necessarily vanishes in
zeroth order according to the commutativity of the undeformed algebra.

In the case when we have a specific momentum map we can simplify the
groups $\Uha$ and $\Glha$ further. In fact, only central elements are
needed. To see this, we proceed in several steps:
\begin{lemma}
    \label{lemma:NicerUhaGlha}%
    Let $\mathcal{A}$ be a unital $^*$-algebra with a $^*$-action of
    $H$ and a momentum map $J$.
    \begin{lemmalist}
    \item A central element $z \in \Zentrum{\mathcal{A}}$ is invariant
        under the action, i.e. $g \acts z = \epsilon(g) z$ for all $g
        \in H$.
    \item The group $\Uhza$ coincides with the unital
        $^*$-homomorphisms $\mathsf{z}: H \longrightarrow
        \Zentrum{\mathcal{A}}$.
    \item The map
        \begin{equation}
            \label{eq:AdJIsomorphism}
            \Uha \ni \mathsf{a}
            \longmapsto 
            \mathsf{a}_J = J^{-1} \ast \mathsf{a} \ast J
            \in \Uhza
        \end{equation}
        is an isomorphism of groups. Here $J^{-1}(g) = J(S(g))$ is the
        convolution inverse of $J$.
    \end{lemmalist}
    In the ring-theoretic situation the analogous statement holds with
    $\Uha$ and $\Uhza$ being replaced by $\Glha$ and $\Glhza$,
    respectively.
\end{lemma}
\begin{proof}
    The first part was already noted in
    \eqref{eq:CentralElementsInvariant}. Since $\Zentrum{\mathcal{A}}$
    is commutative by definition, we can apply
    Lemma~\ref{lemma:TrivialActionOnCommutativeAlgebra} to obtain the
    second part.  For the last part, we show that $\mathsf{a}_J(g)$ is
    central for every $g \in H$. Indeed, for all $b \in \mathcal{A}$
    we have
    \begin{align*}
        \mathsf{a}_J(g) b
        &=
        J(S(g_\sweedler{1})) \mathsf{a}(g_\sweedler{2})
        J(g_{\sweedler{3}}) b \\
        &=
        J(S(g_\sweedler{1})) \mathsf{a}(g_\sweedler{2})
        (g_\sweedler{3} \acts b) J(g_\sweedler{4}) \\
        &=
        J(S(g_\sweedler{1})) (g_\sweedler{2} \acts b)
        \mathsf{a}(g_\sweedler{3}) J(g_\sweedler{4}) \\
        &=
        b J(S(g_\sweedler{1})) \mathsf{a}(g_\sweedler{2})
        J(g_\sweedler{3}) \\
        &=
        b \mathsf{a}_J(g),
    \end{align*}
    using the fact that we have an inner action and using
    Definition~\ref{definition:GroupsGlHAUHA},
    \refitem{item:Commutator}, for $\mathsf{a}$. Moreover,
    $\mathsf{a}_J(\Unit_H) = \Unit_{\mathcal{A}}$ is clear. Since
    $\mathsf{a}_J$ takes values in the center only, one easily
    computes $\mathsf{a}_J(gh) = \mathsf{a}_J(g) \mathsf{a}_J(h)$
    using the property \refitem{item:Product} from
    Definition~\ref{definition:GroupsGlHAUHA}. It follows that
    $\mathsf{a}_J$ is a unital homomorphism into
    $\Zentrum{\mathcal{A}}$.  Note that $J$ has a convolution inverse
    given by $J^{-1}(g) = J(S(g))$ since it is a unital
    homomorphism. From this and \eqref{eq:ConvolutionInvolution} we
    see that $J^* = J^{-1}$. It easily follows that
    \eqref{eq:AdJIsomorphism} maps convolution products to convolution
    products and hence gives a group morphism into $\Glhza$. Finally,
    $(\mathsf{a}^*)_J = (\mathsf{a}_J)^*$ holds which shows that
    $\mathsf{a}_J \in \Uhza$ for $\mathsf{a} \in \Uha$. The
    injectivity of \eqref{eq:AdJIsomorphism} is clear, we only have to
    check surjectivity. Thus let $\mathsf{z} \in \Uhza$ be given and
    consider the pre-image $\mathsf{a} = J * \mathsf{z} * J^{-1}$
    under \eqref{eq:AdJIsomorphism} in $\Hom_{\ring{C}}(H,
    \mathcal{A})$. Clearly, $\mathsf{a}(\Unit_H) =
    \Unit_{\mathcal{A}}$ holds. Using the facts that $\mathsf{z}$
    takes values in the center and that the action is inner allows to
    show the properties \refitem{item:Product} and
    \refitem{item:Commutator} for $\mathsf{a}$ by a simple
    computation. Thus $\mathsf{a} \in \Glha$ follows. Finally, from
    $\mathsf{z}^{-1} = \mathsf{z}^*$ and the compatibility of
    \eqref{eq:AdJIsomorphism} with the $^*$-involution
    \eqref{eq:ConvolutionInvolution} we conclude that $\mathsf{a} \in
    \Uha$ as wanted. This concludes the proof as the ring-theoretic
    part was done on the way as well.
\end{proof}
The next lemma characterizes the uniqueness properties of momentum
maps which lead to the same action:
\begin{lemma}
    \label{lemma:UniquenessMomentumMap}%
    Let $\mathcal{A}$ be a unital $^*$-algebra with $^*$-action of a
    Hopf $^*$-algebra $H$.  For any two momentum maps $J, J' \in
    \starMM_H(\mathcal{A})$ the map $\mathsf{z} = J^{-1} * J'$ defines
    an element $\mathsf{z} \in \Uhza$. Conversely, for $J \in
    \starMM_H(\mathcal{A})$ and $\mathsf{z} \in \Uhza$ we have $J \ast
    \mathsf{z} \in \starMM_H(\mathcal{A})$.
\end{lemma}
\begin{proof}
    Let $a \in \mathcal{A}$ then we compute
    \[
    \mathsf{z}(g)a
    =
    J(S(g_\sweedler{1})) J'(g_\sweedler{2})a
    =
    J(S(g_\sweedler{1})) (g_\sweedler{2} \acts a) J'(g_\sweedler{3})
    =
    a J(S(g_\sweedler{1})) J'(g_\sweedler{2})
    =
    a\mathsf{z}(g)
    \]
    showing $\mathsf{z}(g)\in \Zentrum{\alg{A}}$ for all $g \in H$.
    Using \eqref{eq:CentralElementsInvariant}, the remaining
    properties for $\mathsf{z} \in \Uhza$ are easily verified and the
    converse statement follows directly.
\end{proof}
The ring-theoretic situation is handled analogously replacing $\Uhza$
by $\Glhza$: we do not need to formulate this in detail.  In the case
of $(C^{\infty}(M)[[\lambda]], \star)$ and
$\mathcal{U}_{\lambda}(\lie{g})$, the ambiguity in the existence of
quantum momentum maps reduces to the one in
\cite{mueller-bahns.neumaier:2004a}, see also the discussion in
\cite{waldmann:2006a}.

Since $\Uhza$ and $\Uha$ are isomorphic via \eqref{eq:AdJIsomorphism}
we can also parametrize momentum maps $J'$ starting with a fixed
reference momentum map $J$ and $\mathsf{a} \in \Uha$ by $J' = J *
\mathsf{z} = \mathsf{a} * J$ where $\mathsf{z} = \mathsf{a}_J \in
\Uhza$. We summarize these results in the following theorem:
\begin{theorem}
    \label{theorem:HowManyMomentumMaps}%
    Let $\mathcal{A}$ be a unital $^*$-algebra with $^*$-action of $H$
    and a momentum map $J$. Then the group $\Uha$ acts freely and
    transitively from the left on $\starMM_H(\mathcal{A})$ by left
    multiplication using the convolution product $*$ of
    $\Hom_{\ring{C}}(H, \mathcal{A})$. Equivalently, $\Uhza$ acts
    freely and transitively on $\starMM_H(\mathcal{A})$ from the right
    by the convolution product. In the ring-theoretic situation,
    $\starMM_H(\mathcal{A})$, $\Uha$, and $\Uhza$ have to be replaced
    by $\MM_H(\mathcal{A})$, $\Glha$, and $\Glhza$, respectively.
\end{theorem}

%
%

\subsection{The splitting}
\label{subsec:Splitting}

The following consideration is central. We state the result in
slightly larger generality though we need this only for equivalence
bimodules later.
\begin{proposition}
    \label{proposition:impulsabbildungenLiftBimodul}%
    Consider $^*$-algebras $\mathcal{A}$, $\mathcal{B}$, and
    $\mathcal{C}$ with $^*$-actions by $H$ and momentum maps
    $J^{\mathcal{A}}$, $J^{\mathcal{B}}$, and $J^{\mathcal{C}}$,
    respectively.
    \begin{propositionlist}
    \item If $\BEA$ is a $(\mathcal{B}, \mathcal{A})$-bimodule then
        \begin{equation}
            \label{eq:LiftedAction}
            h \acts x
            =
            J^{\mathcal{B}}(h_\sweedler{1})
            \cdot x \cdot
            J^{\mathcal{A}}(S(h_\sweedler{2}))
        \end{equation}
        for $h \in H$ and $x \in \BEA$ defines an action of $H$ on
        $\BEA$ which is compatible with the bimodule structure.
    \item Endowing two $(\mathcal{B}, \mathcal{A})$-bimodules $\BEA$
        and $\BEpA$ with this $H$-action makes every bimodule morphism
        $\Phi: \BEA \longrightarrow \BEpA$ equivariant with respect to
        the $H$-actions \eqref{eq:LiftedAction}.
    \item If a $(\mathcal{B}, \mathcal{A})$-bimodule $\BEA$ carries an
        algebra-valued inner product then the action
        \eqref{eq:LiftedAction} is compatible with it.
    \item For a $(\mathcal{C}, \mathcal{B})$-bimodule $\CFB$ and a
        $(\mathcal{B}, \mathcal{A})$-bimodule $\BEA$ the lift
        \eqref{eq:LiftedAction} to $\CFB \tensor[\mathcal{B}] \BEA$
        coincides with the tensor product of the lifts to $\CFB$ and
        $\BEA$.
    \item For the canonical bimodule $\AAA$ the lift
        \eqref{eq:LiftedAction} reproduces the original action of $H$
        on $\mathcal{A}$.
    \end{propositionlist}
\end{proposition}
\begin{proof}
    It follows immediately that \eqref{eq:LiftedAction} yields a
    compatible $H$-action on the bimodule $\BEA$. Being an ``inner''
    action the second statement is trivial.  The compatibility with
     either a $\mathcal{B}$-valued or an $\mathcal{A}$-valued inner
    product is a direct calculation
    \begin{align*}
        h \acts \SPA{x,y}
        &=
        J^{\mathcal{A}}(h\sweedler{1}) \SPA{x,y}
        J^{\mathcal{A}}(S(h \sweedler{2}))\\
        &=
        J^{\mathcal{A}}(h\sweedler{1})
        \SPA{
          J^{\mathcal{B}}(S(h\sweedler{2})^{*}) \cdot
          x,J^{\mathcal{B}}(h\sweedler{3}) \cdot y
        }
        J^{\mathcal{A}}(S(h \sweedler{4}))\\
        &=
        \SPA{
          J^{\mathcal{B}}(S(h\sweedler{2})^{*}) \cdot x \cdot
          J^{\mathcal{A}}(S(S(h\sweedler{1})^{*})),
          J^{\mathcal{B}}(h\sweedler{3}) \cdot y \cdot
          J^{\mathcal{A}}(S(h\sweedler{4}))
        }\\
        &=
        \SPA{S(h\sweedler{1})^{*}\acts x,h\sweedler{2}\acts y}.
    \end{align*}
    The compatibility with the tensor product follows from
    \begin{align*}
        (h_\sweedler{1} \acts x)
        \tensor[\alg{B}]
        (h_\sweedler{2} \acts y)
        &=
        \left(
            J^{\alg{C}}(h_\sweedler{1})
            \cdot x \cdot
            J^{\alg{B}}(S(h_\sweedler{2}))
        \right)
        \tensor[\alg{B}]
        \left(
            J^{\alg{B}}(h_\sweedler{3})
            \cdot y \cdot
            J^{\alg{A}}(S(h_\sweedler{4}))
        \right)\\
        &=
        \left(J^{\alg{C}}(h_\sweedler{1}) \cdot x\right)
        \tensor[\alg{B}]
        \left(
            y \cdot J^{\alg{A}}(h_\sweedler{2})
        \right)\\
        &=
        h \acts (x \tensor[\alg{B}]{y}).
    \end{align*}
    The last part is clear.
\end{proof}
Note that the proposition also holds in the ring-theoretic case,
except of course for the statement about the inner products.

In the following, we consider $H$-equivariant Morita theory only for
$^*$-algebras \emph{with} a specific momentum map. In particular, we
consider the Picard groupoids of the various flavours only over this
restricted class of algebras. From the above
Proposition~\ref{proposition:impulsabbildungenLiftBimodul} it follows
that for every $^*$- or strong equivalence bimodule $\BEA$ the lift
\eqref{eq:LiftedAction} endows $\BEA$ with the structure of an
$H$-equivariant $^*$- or strong equivalence bimodule,
respectively. Moreover, the compatibility with tensor products and
morphisms guarantees that we obtain a \emph{groupoid morphism} after
passing to isomorphism classes. We denote these groupoid morphisms by
\begin{equation}
    \label{eq:DefinitionJforstarStr}
    J: \starPic \longrightarrow \starPicH
    \quad
    \textrm{and}
    \quad
    J: \StrPic \longrightarrow \StrPicH,
\end{equation}
respectively. Note that we suppress the restriction to this more
specific class of $^*$-algebras in our notation. Analogously, in the
ring-theoretic situation we get a groupoid morphism
\begin{equation}
    \label{eq:JforRingTheoreticPic}
    J: \Pic \longrightarrow \PicH.
\end{equation}
\begin{remark}
    \label{remark:JBimodToBimodH}%
    In fact, the above proposition even gives a functor $J: \BiMod
    \longrightarrow \BiModH$ whose restriction to invertible arrows is
    \eqref{eq:JforRingTheoreticPic} and similarly for the $^*$- and
    strong version. However, we will not need this additional
    structure.
\end{remark}
The following theorem is now an easy consequence of these
considerations which we formulate for all three versions of the Picard
groupoid:
\begin{theorem}
    \label{theorem:splitting the image starpich to starpic}%
    Consider the Picard groupoids over the restricted class of unital
    $^*$-algebras (or algebras) with momentum maps.
    \begin{theoremlist}
    \item \label{item:SurjectiveGroupoidMorphism} The groupoid
        morphisms $J$ from \eqref{eq:DefinitionJforstarStr} and
        \eqref{eq:JforRingTheoreticPic}, respectively, are right
        inverses of the canonical groupoid morphisms in
        \eqref{eq:NaturalForgetting}. In particular, over this class,
        \eqref{eq:NaturalForgetting} is surjective.
    \item \label{item:QuotientUnnecessary} Let $\mathcal{A}$ be fixed
        unital $^*$-algebra (or algebra) with momentum map. The images
        of the subgroups $\Glza$ and $\Uza$ in $\Glha$ and $\Uha$,
        respectively, under the morphism \eqref{eq:von GlZ(A) nach
          GL(H,A)} are trivial. Hence we identify $\Gloha=\Glha$ and
        $\Uoha=\Uha$.
    \item \label{item:SplitExactSequences} For a fixed unital
        $^*$-algebra (or algebra) with momentum map there are split
        exact sequences
        \begin{equation}
            \label{eq:Spaltende  stern-Bimod exakte sequenz}
            \xy
            \morphism(0,0)<700,0>[1`\Uha;]
            \morphism(700,0)<900,0>[\Uha`\starPicH(\alg{A});]
            \morphism(1600,0)<1000,0>[\starPicH(\alg{A})`\starPic(\alg{A});]
            \morphism(2600,0)<800,0>[\starPic(\alg{A})`1;]
            \morphism(1600,0)|a|/{@{<-}@/_2em/}/<1000,0>[\starPicH(\alg{A})`\starPic(\alg{A});J]
            \endxy,
        \end{equation}
        \begin{equation}
            \label{eq:SpaltendeStr-Bimod exakte sequenz}
            \xy
            \morphism(0,0)<700,0>[1`\Uha;]
            \morphism(700,0)<900,0>[\Uha`\StrPicH(\alg{A});]
            \morphism(1600,0)<1000,0>[\StrPicH(\alg{A})`\StrPic(\alg{A});]
            \morphism(2600,0)<800,0>[\StrPic(\alg{A})`1;]
            \morphism(1600,0)|a|/{@{<-}@/_2em/}/<1000,0>[\StrPicH(\alg{A})`\StrPic(\alg{A});J]
            \endxy,
        \end{equation}
        and
        \begin{equation}
            \label{eq:Spaltende  Bimod exakte sequenz}
            \xy
            \morphism(0,0)<700,0>[1`\Glha;]
            \morphism(700,0)<900,0>[\Glha`\PicH(\alg{A});]
            \morphism(1600,0)<1000,0>[\PicH(\alg{A})`\Pic(\alg{A});]
            \morphism(2600,0)<800,0>[\Pic(\alg{A})`1;]
            \morphism(1600,0)|a|/{@{<-}@/_2em/}/<1000,0>[\PicH(\alg{A})`\Pic(\alg{A});J]
            \endxy,
        \end{equation}    
        respectively.
    \item \label{item:Semidirect} For a fixed unital $^*$-algebra (or
        algebra) with momentum map the $H$-equivariant Picard groups
        are semidirect products
        \begin{equation}
            \label{eq:SemidirectStarVersion}
            \starPicH(\alg{A}) \cong \starPic(\alg{A}) \ltimes \Uha
            \quad
            \textrm{and}
            \quad
            \StrPicH(\alg{A}) \cong \StrPic(\alg{A}) \ltimes \Uha,
        \end{equation}
        as well as
        \begin{equation}
            \label{eq:SemidirectRingVersion}
            \PicH(\alg{A}) \cong \Pic(\alg{A}) \ltimes \Glha,
        \end{equation}
        respectively. The group actions responsible for the semidirect
        product are induced by 
        \begin{equation}
            \label{eq:SemidirectAction}
            \Phi(\mathcal{E}): \Uha
            \ni \mathsf{a} 
            \longmapsto
            h_{J(\mathcal{E})} (\mathsf{a}) \in
            \Uha
        \end{equation}
        for the $^*$- and strong version with $h_{J(\mathcal{E})}$ as
        in Lemma~ \ref{lemma:Isomorphism uha uhb}, and by
        \begin{equation}
            \label{eq:SemidirectActionRing}
            \Phi(\mathcal{E}): \Glha
            \ni \mathsf{a} 
            \longmapsto
            h_{J(\mathcal{E})} (\mathsf{a}) \in
            \Glha
        \end{equation}
        in the ring-theoretic framework, respectively.
    \end{theoremlist}
\end{theorem}
\begin{proof}
    The first part is clear as forgetting the lifted action reproduces
    the bimodule we started
    with. Proposition~\ref{proposition:impulsabbildungenLiftBimodul}
    ensures that $J$ is indeed a groupoid morphism. The second part
    follows from $ \Zentrum{\mathcal{A}} \subseteq \mathcal{A}^H$ for
    an algebra with an inner $H$-action. The third part follows from
    \eqref{eq:kernel picH-pic} and \eqref{eq:kernel starpicH-starpic},
    respectively, and the groupoid structure at once. For the last
    part we have to compute the precise action of the Picard group on
    the groups $\Uha$ and $\Glha$, respectively.  Thus let $\AEA$ be a
    $^*$-equivalence $(\mathcal{A}, \mathcal{A})$-bimodule and
    $\mathsf{a} \in \Uha$. In order to compute the action of $\AEA$ on
    $\mathsf{a}$ we have to consider the $^*$-equivalence bimodule
    $\AEA \tensor \AAA \tensor \AccEA \cong \AAA$. When $\AEA$ is
    equipped with the lifted action we have to determine the action
    $\acts^{\Phi(\mathcal{E})(\mathsf{a})}$ on this tensor product
    when $\AAA$ was endowed with the twisted action
    $\acts^{\mathsf{a}}$. But this is now an easy computation as the
    action $\acts^{\Phi(\mathcal{E})(\mathsf{a})}$ is given by the
    tensor product of the actions on each factor, respectively.  Let
    $x \otimes a \otimes \cc{y} \in J(\mathcal{E}) \otimes \alg{A}
    \otimes J(\cc{\mathcal{E}})$. Then
   \begin{align*}
       g \acts^{\Phi(\mathcal{E}) (\mathsf{a})} 
       \left(x \otimes a \otimes \cc{y}\right)
       &=
       (g_\sweedler{1} \acts x) 
       \otimes
       (g_\sweedler{2} \acts^{\mathsf{a}} a)
       \otimes
       (g_\sweedler{3} \ccacts \cc{y}) \\
       &=
       (g_\sweedler{1} \acts x)
       \otimes
       \left(
           \mathsf{a}(g_\sweedler{2}) \cdot (g\sweedler{3} \acts a)
       \right)
       \otimes
       (g_\sweedler{4} \ccacts \cc{y}) \\
       &=
       \left(
           (g_\sweedler{1} \acts x) \cdot \mathsf{a}(g_\sweedler{2})
       \right)
       \otimes
       (g\sweedler{3} \acts a)
       \otimes
       (g_\sweedler{4} \ccacts \cc{y}) \\
       &=
       (g_\sweedler{1} \acts_{\mathsf{a}} x)
       \otimes
       (g_\sweedler{2} \acts a)
       \otimes
       (g_\sweedler{3} \ccacts \cc{y}) \\
       &=
       (g_\sweedler{1} \acts^{h_{J(\mathcal{E})}(\mathsf{a})} x)
       \otimes
       (g_\sweedler{2} \acts a)
       \otimes
       (g_\sweedler{3} \ccacts \cc{y}) \\
       &=
       g \acts^{h_{J(\mathcal{E})}(\mathsf{a})}
       (x \otimes a \otimes \cc{y}).
   \end{align*}
   Thus, $\Phi(\mathcal{E}) (\mathsf{a}) =
   h_{J(\mathcal{E})}(\mathsf{a})$ follows. The ring-theoretic case is
   analogous.
\end{proof}

On the other hand, the existence of a momentum map is an invariant for
equivariant Morita equivalence. This shows that our restriction to the
above class is not severe: all other unital $^*$-algebras (algebras)
with other types of $^*$-actions (actions) of $H$ are necessarily on
different connected components of the corresponding Picard groupoids
and hence ``invisible'' from the Morita theory point of view. In
\cite{jansen.waldmann:2006a}, Morita invariants are seen as arising
from actions of the corresponding Picard groupoids. We follow this
point of view and formulate this Morita invariant by means of an
action. To this end we need the following technical lemma:
\begin{lemma}
    \label{lemma:Existence of momentum map is a PicH invariant}
    Consider unital $^*$-algebras $\alg{A}$, $\alg{B}$, and
    $\mathcal{C}$ with a $^*$-action of the Hopf $^*$-algebra
    $H$. Furthermore, assume that $\starMM_H(\mathcal{A}) \ne
    \emptyset$.
    \begin{lemmalist}
    \item Assume there is an $H$-equivariant $^*$-equivalence bimodule
        $\BEA$. For $J \in \starMM_H(\mathcal{A})$ there exists a
        uniquely determined momentum map $h_{\mathcal{E}}(J) \in
        \starMM_H(\mathcal{B})$ such that for all $g \in H$ and $x \in
        \BEA$
        \begin{equation}
            \label{eq:InducedMomentumMapForB}
            g \acts x
            =
            h_{\mathcal{E}}(J)(g_\sweedler{1})
            \cdot x \cdot
            J(S(g_{2})).
        \end{equation}
    \item Assume $\BEA, \BEpA$ are isomorphic $H$-equivariant
        $^*$-equivalence bimodules. Then
        \begin{equation}
            \label{eq:hEhEpEqual}
            h_{\mathcal{E}} = h_{\mathcal{E}'}.
        \end{equation}
    \item Assume $\BEA$ and $\CFB$ are $H$-equivariant
        $^*$-equivalence bimodules. Then
        \begin{equation}
            \label{eq:hFhEIshFEhAIdentity}
            h_{\mathcal{F}} \circ h_{\mathcal{E}}
            = h_{\mathcal{F} \tensor \mathcal{E}}
            \quad
            \textrm{and}
            \quad
            h_{\mathcal{A}} = \id_{\starMM_H(\mathcal{A})}.
        \end{equation}
    \item Assume $\BEA$ is an $H$-equivariant $^*$-equivalence
        bimodule. Then for $\mathsf{a} \in \Uha$ and $J \in
        \starMM_H(\mathcal{A})$ we have
        \begin{equation}
            \label{eq:hECompatibleWithAllStuff}
            h_{\mathcal{E}}(\mathsf{a} * J) =
            h_{\mathcal{E}}(\mathsf{a}) * h_{\mathcal{E}}(J).
        \end{equation}
    \end{lemmalist}
\end{lemma}
\begin{proof}
    For the first part, let $a \in \mathcal{A}$ and $x \in \BEA$. Then
    for $g \in H$ we have
    \begin{align*}
        (g_\sweedler{1} \acts (x \cdot a) )\cdot J(g_\sweedler{2})
        &=
        (g_\sweedler{1} \acts x)
        \cdot \left(
            (g_\sweedler{2} \acts a) J(g_\sweedler{3})
        \right) \\
        &=
        (g_\sweedler{1} \acts x)
        \cdot \left(
            J(g_\sweedler{2}) a J(S(g_\sweedler{3})) J(g_\sweedler{4})
        \right) \\
        &=
        \left(
            (h_\sweedler{1} \acts x ) \cdot J(h_\sweedler{2})
        \right)
        \cdot a.
    \end{align*}
    This shows that the map $x \mapsto (g_\sweedler{1} \acts x) \cdot
    J(g_\sweedler{2})$ is right $\mathcal{A}$-linear and depends
    linearly on $g$. Thus it is the left multiplication by a unique
    element $h_{\mathcal{E}}(J)(g) \in \mathcal{B}$, which depends
    linearly on $g$. We claim that this defines the momentum map we
    are looking for. First we have for $g, h \in H$ and $x \in \BEA$
    \[
    \left(
        h_{\mathcal{E}}(J)(g) h_{\mathcal{E}}(J)(h)
    \right) \cdot x 
    =
    h_{\mathcal{E}}(J)(gh) \cdot x,
    \]
    by simply expanding the definition. Moreover,
    $h_{\mathcal{E}}(J)(\Unit_H) = \Unit_{\mathcal{B}}$ is clear. It
    is now a simple check that \eqref{eq:InducedMomentumMapForB}
    holds. Let $b \in \mathcal{B}$ and $x \in \BEA$ be given. Then
    \begin{align*}
        (g \acts b) \cdot x
        &=
        g_\sweedler{1} \acts
        \left(b \cdot (S(g_\sweedler{2}) \acts x)\right) \\
        &=
        h_{\mathcal{E}}(J)(g_\sweedler{1}) \cdot
        \left(
            b \cdot
            \left(
                h_{\mathcal{E}}(J)(S(g_\sweedler{3})_\sweedler{1})
                \cdot x \cdot
                J(S(S(g_\sweedler{3})_\sweedler{2}))
            \right)
        \right)
        \cdot
        J(S(g_\sweedler{2})) \\
        &=
        \left(
            h_{\mathcal{E}}(J)(g_\sweedler{1})
            b
            h_{\mathcal{E}}(J)(S(g_\sweedler{4}))
        \right)
        \cdot x \cdot
        \left(
            J(S^2(g_\sweedler{3})) J(S(g_\sweedler{2}))
        \right) \\
        &=
        \left(
            h_{\mathcal{E}}(J)(g_\sweedler{1})
            b
            h_{\mathcal{E}}(J)(S(g_\sweedler{4}))
        \right)
        \cdot x \cdot
        \left(
            J(S(g_\sweedler{2} S(g_\sweedler{3})))
        \right) \\
        &=
        \left(
            h_{\mathcal{E}}(J)(g_\sweedler{1})
            b
            h_{\mathcal{E}}(J)(S(g_\sweedler{2}))
        \right)
        \cdot x
    \end{align*}
    shows that $h_{\mathcal{E}}(J)$ is a momentum map for the algebra
    $\mathcal{B}$. Finally, we use the compatibility of the
    $\mathcal{A}$-valued inner product on $\BEA$ with the action and
    compute for $x, y \in \BEA$ and $g \in H$
    \begin{align*}
        \SPA{h_{\mathcal{E}}(J)(g)^* \cdot x, y}
        &=
        \SPA{x, h_{\mathcal{E}}(J) \cdot y} \\
        &=
        \SPA{x, (g_\sweedler{1} \acts y) \cdot J(g_\sweedler{2})} \\
        &=
        \left(
            g_\sweedler{2} \acts \SPA{g_\sweedler{1}^* \acts x, y}
        \right)
        J(g_\sweedler{3}) \\
        &=
        J(g_\sweedler{2}) \SPA{g_\sweedler{1}^* \acts x , y} \\
        &=
        \SPA{
          (g_\sweedler{1}^* \acts x) \cdot J(g_{\sweedler{2}}^*),
          y
        } \\
        &=
        \SPA{h_{\mathcal{E}}(J)(g^*) \cdot x, y}.
    \end{align*}
    Since the inner product is non-degenerate this establishes
    $h_{\mathcal{E}}(J)(g)^* = h_{\mathcal{E}}(J)(g^*)$. Thus
    $h_{\mathcal{E}}(J)$ is indeed a momentum map for the $^*$-algebra
    $\mathcal{B}$. Finally, the way we constructed
    $h_{\mathcal{E}}(J)$ gives the uniqueness. This completes the
    first part. For the second part, let $\Psi: \BEA \longrightarrow
    \BEpA$ be an $H$-equivariant isometric isomorphism. Then
    \[
    g \acts x
    =
    \Psi^{-1}\left(
        g \acts \Psi(x)
    \right)
    =
    \Psi^{-1}\left(
        h_{\mathcal{E}'}(J)(g_\sweedler{1})
        \cdot \Psi(x) \cdot
        J(S(g_\sweedler{2}))
    \right)
    =
    h_{\mathcal{E}'}(J)(g_\sweedler{1})
    \cdot x \cdot
    J(S(g_\sweedler{2}))
    \]
    by equivariance and the fact that $\Psi$ is a bimodule
    isomorphism. This shows that $h_{\mathcal{E}'}(J)$ satisfies
    \eqref{eq:InducedMomentumMapForB} and thus coincides with
    $h_{\mathcal{E}}(J)$ by uniqueness. For the third part, let $\CFB$
    be yet another $H$-equivariant $^*$-equivalence bimodule. For
    factorizing tensors $y \tensor x \in \CFB \tensor \BEA$ we unwind
    the definition
    \begin{align*}
        h_{\mathcal{F} \tensor \mathcal{E}}(J) \cdot (y \tensor x)
        &=
        (g_\sweedler{1} \acts (y \tensor x)) \cdot J(g_\sweedler{2})
        \\
        &=
        (g_\sweedler{1} \acts y) \tensor
        \left(
            (g_\sweedler{2} \acts x) \cdot J(g_\sweedler{3})
        \right) \\
        &=
        (g_\sweedler{1} \acts y) \tensor
        \left(
            h_{\mathcal{E}}(J)(g_\sweedler{2}) \cdot x
        \right) \\
        &=
        h_{\mathcal{F}}(h_{\mathcal{E}}(J))(g) \cdot (y \tensor x),
    \end{align*}
    showing the first claim in \eqref{eq:hFhEIshFEhAIdentity}. The
    second statement in \eqref{eq:hFhEIshFEhAIdentity} is trivial. For
    the last part we consider $\mathsf{a} \in \Uha$ and compute for $x
    \in \BEA$
    \begin{align*}
        h_{\mathcal{E}}(\mathsf{a} * J) \cdot x
        &=
        (g_\sweedler{1} \acts x) \cdot \left(
            \mathsf{a}(g_\sweedler{2}) J(g_\sweedler{3})
        \right) \\
        &=
        (g_\sweedler{1} \acts_{\mathsf{a}} x) \cdot J(g_\sweedler{2}) \\
        &=
        (g_\sweedler{1} \acts^{h_{\mathcal{E}}(\mathsf{a})} x)
        \cdot J(g_\sweedler{2}) \\
        &=
        h_{\mathcal{E}}(\mathsf{a})(g_\sweedler{1}) \cdot
        \left(
            (g_\sweedler{2} \acts x) \cdot J(g_\sweedler{3})
        \right) \\
        &=
        \left(
            h_{\mathcal{E}}(\mathsf{a})(g_\sweedler{1})
            h_{\mathcal{E}}(J)(g_\sweedler{2})
        \right)
        \cdot x,
    \end{align*}
    using Lemma~\ref{lemma:Isomorphism uha uhb}. This shows
    \eqref{eq:hECompatibleWithAllStuff}.
\end{proof}

Rephrasing the above lemma gives us the following general statement on
the existence of momentum maps:
\begin{theorem}[$\starPicH$ acts on $\starMM_H$]
    \label{theorem:PicHActsOnMoMa}%
    The $H$-equivariant $^*$-Picard groupoid $\starPicH$ acts on the
    $\mathsf{U}(H, \cdot)$-spaces $\starMM_H(\cdot)$ from the
    left. In particular, the existence of a momentum map is invariant
    under $H$-equivariant $^*$-Morita equivalence.
\end{theorem}
Here it is understood that the groupoid action of $\starPicH$ on
$\starMM_H(\cdot)$ is ``along'' the groupoid action of $\starPicH$
on $\mathsf{U}(H, \cdot)$ from Lemma~\ref{lemma:Isomorphism uha
  uhb}. In formulas, this just means
\eqref{eq:hECompatibleWithAllStuff}.

Dropping the $^*$-involutions gives us the ring-theoretic result which
means that the $H$-equivariant Picard groupoid acts on the
$\mathsf{Gl}(H, \cdot)$-spaces $\MM_H(\cdot)$ from the left. Again,
the existence of a momentum map is invariant under $H$-equivariant
Morita equivalence.

In the sequel, we will need more precisely the structure of
$\starPicH(\alg{A})$ for a \emph{commutative} unital $^\ast$-algebra
$\alg{A}$. In this case, we do, of course, not assume to have an inner
action.
\begin{lemma}
    \label{lemma:HEonCenter}%
    Let $\alg{A}$ and $\alg{B}$ be commutative unital $^\ast$-algebras
    and let $[\BEA] \in \starPicH(\alg{B},\alg{A})$. Then there is an
    $H$-equivariant $^\ast$-isomorphism $h_{\mathcal{E}}: \alg{A}
    \longrightarrow \alg{B}$ determined by $h_{\mathcal{E}}(a) \cdot x
    = x \cdot a$. The map $\mathcal{E} \mapsto h_{\mathcal{E}}$ is
    compatible with the composition and units in $\starPicH$.
\end{lemma}
\begin{proof}
    This is \cite [Lemma 5.1]{jansen.waldmann:2006a}, except for the
    $H$-equivariance of $h_{\mathcal{E}}$ which is easy to see.
\end{proof}

Denote by $\starAutH(\alg{A})$ the $H$-equivariant
$^\ast$-automorphism of the $^\ast$-algebra $\alg{A}$ and by
$\starSPicH(\alg{A})$ the equivalence classes of static
$H$-equivariant bimodules in $\starPicH(\alg{A})$. Recall that a
$(\alg{A}, \alg{A})$-bimodule $\AEA$ is called \emph{static} (or
commutative or central), if $a \cdot x = x \cdot a$ for all $a \in
\alg{A}$ and $x \in \AEA$. Analogously, one defines
$\StrSPicH(\mathcal{A})$ and the ring-theoretic version
$\SPicH(\mathcal{A})$.
\begin{proposition}
    \label{proposition:semidirect product for starPicH}%
    For a commutative unital $^\ast$-algebra $\alg{A}$ with
    $^*$-action of $H$ one has
    \begin{equation}
        \label{eq:semidirect product for starPicH}
        \starPicH(\alg{A})
        \cong
        \starAutH(\alg{A}) \ltimes \starSPicH(\alg{A}),
    \end{equation}
    and analogously for $\StrPicH(\alg{A})$. In the ring-theoretic
    case we have $\PicH(\mathcal{A}) \cong \AutH(\mathcal{A}) \ltimes
    \SPicH(\mathcal{A})$.
\end{proposition}
\begin{proof}
    For an $^*$-automorphism $\Phi \in \starAut(\mathcal{A})$ we can
    twist the canonical strong self-equivalence $\AAA$ by modifying
    the left module structure to $a \cdot_\Phi x = \Phi^{-1}(a)
    x$. The right module structure is kept unchanged and the inner
    products are changed appropriately. This gives a group morphism
    $\ell: \starAut(\mathcal{A}) \longrightarrow
    \StrPic(\mathcal{A})$. It turns out that this is compatible with
    the $H$-equivariance if $\Phi$ is $H$-equivariant. In fact, this
    construction comes from a groupoid morphism $\ell$ from the
    $H$-equivariant $^*$-isomorphism groupoid $\starIsoH$ to
    $\StrPicH$, see \cite[Prop.~4.7]{jansen.waldmann:2006a}.  If now
    $\mathcal{A}$ is commutative $\StrSPicH(\mathcal{A})$ can be
    viewed as a normal subgroup of $\StrPicH(\mathcal{A})$ and with
    the map $h$ from Lemma~\ref{lemma:HEonCenter} we get an exact
    sequence
    \[
    \xy
    \morphism(0,0)<700,0>[1`\StrSPicH(\mathcal{A});]
    \morphism(700,0)<900,0>[\StrSPicH(\mathcal{A})`\StrPicH(\alg{A});]
    \morphism(1600,0)<1000,0>[\StrPicH(\alg{A})`\starAutH(\alg{A});h]
    \morphism(2600,0)<800,0>[\starAutH(\alg{A})`1;]
    \morphism(1600,0)|a|/{@{<-}@/_2em/}/<1000,0>[\StrPicH(\alg{A})`\starAutH(\alg{A});\ell]
    \endxy,
    \]
    which splits by $\ell$. The case of $^*$-equivalence is analogous
    and the ring-theoretic case is even simpler. Note that the above
    split exact sequence is classic for the non-equivariant,
    ring-theoretic case, see e.g. \cite[Chap.~II,
    Prop.~5.4]{bass:1968a}.
\end{proof}

%
%

\section{Equivariant Morita equivalence of star products}
\label{sec:EMEStarProducts}

We will now apply the notion of equivariant Morita equivalence of
algebras to star products on $(M, \omega)$. The setting is again as in
Section~\ref{subsec:FedosovStarProducts} and we will call the
$\mathcal{U}_{\mathbb{C}}(\lie{g})[[\lambda]]$-equivariance simply
$\lie{g}$-equivariance throughout this section.

%
%

\subsection{Classical and semiclassical limits of a
  $\lie{g}$-equivariant equivalence bimodule}
\label{subsec:semiclass limit}

In the sequel we want to explore the classical and semiclassical limit
of equivalence bimodules with respect to the symmetry structures. We
have to take care of the ``up to isomorphism'' statements as in
Section~\ref{subsec:FedosovStarProducts} more carefully.  Recall that
the classical limit of an equivalence bimodule $\deform{\mathcal{L}}$
for two star products is defined by the $\Fkt$-bimodule
\begin{equation}
    \label{eq:classicalLimitBimodul}
    \cl(\deform{\mathcal{L}})
    =
    \deform{\mathcal{L}}/ \lambda \deform{\mathcal{L}},
\end{equation}  
and yields an equivalence bimodule for $\Fkt$. The quotient map is
also referred to as the \emph{classical limit} and will be denoted by
$\cl$. Furthermore, one has the following result, see
e.g. \cite{bursztyn.waldmann:2004a}:
\begin{proposition}
    \label{proposition:clLimitEquivalenceBimoduL}%
    Let $\deform{\mathcal{L}}$ be an equivalence bimodule for the two
    star products $\star$ and $\star'$ where $\star'$ acts from the
    left. Then there exists a unique symplectomorphism $\Psi$ such
    that $^{\Psi}\deform{\mathcal{L}}$ with the twisted left module
    structure by $\Psi$ is an equivalence bimodule for
    $\Psi^{*}(\star')$ and $\star$ and
    $\cl(^{\Psi}\deform{\mathcal{L}})$ is a static equivalence
    bimodule for $\Fkt$ and hence isomorphic to the sections
    $\mathcal{L} = \Schnitte(L)$ of a line bundle $L$.
\end{proposition}
Thus we can replace $\star'$ by $\Psi^*(\star')$ and assume that the
equivalence bimodule has a static equivalence bimodule as classical
limit. A further result from the general ring-theoretic deformation
theory is that in this case where $\deform{\mathcal{L}}$ has a static
classical limit $\mathcal{L}$ we have a $\mathbb{C}[[\lambda]]$-linear
isomorphism $\deform{\mathcal{L}} \cong \mathcal{L}[[\lambda]]$ which
induces a bimodule deformation of $\mathcal{L}$ encoded by the
deformed left and right multiplication laws $\bullet'$ and $\bullet$
as in \eqref{eq:Abulletprimes} and \eqref{eq:RightModuleStructure},
respectively. Thus it will be sufficient to assume that
$\deform{\mathcal{L}} = (\Schnitte(L)[[\lambda]], \bullet', \bullet)$
for the time being \emph{and} to plug in the symplectomorphism $\Psi$
later on by hand.

Let $L$ be a line bundle and let $\deform{\mathcal{L}} =
\mathcal{L}[[\lambda]]$ be a deformation of the sections $\mathcal{L}
= \Schnitte(L)$ into an equivalence bimodule for two star products
$\star$ and $\star'$.  Assume furthermore, that the first order terms
$C_1$ and $C_1'$ of $\star$ and $\star'$ as in \eqref{eq:StarProduct}
coincide. For example, Fedosov star products always have the same
first order term given by $\frac{\I}{2}$ times the Poisson bracket.
From \cite[Prop.~4.3]{bursztyn:2002a} we know that the difference
of the first order terms of the deformed module structures gives rise
to a connection on $L$, called the \emph{induced connection}. For
Hamiltonian vector fields $X_f$ with $f \in \Fkt$ it is defined by
\begin{equation}
    \label{eq:contravariantConn}
    \nabla^L_{X_f}(s) = \I R_1'(f, s) - \I R_1(s, f),
\end{equation}
and by $\Fkt$-linear expansion this gives a well-defined connection
since the Hamiltonian vector fields span the tangent spaces at every
point. Note that in the Poisson case we would only obtain a
\emph{contravariant} connection. If we pass to an equivalent bimodule
deformation by an equivalence transformation this connection does not
change, see \cite[Prop.~4.17]{bursztyn.waldmann:2004a}.

If $\deform{\mathcal{L}}$ was endowed with inner products defining a
$^{*}$-equivalence bimodule for the two star products, the classical
limit yields in addition a pseudo Hermitian fiber metric on the line
bundle in the classical limit. Indeed, one defines the classical limit
of the inner product $\deform{h}$ by
\begin{equation}
    \label{eq:clhDef}
    h(s, s') = \deform{h}(\deform{s}, \deform{s'})\big|_{\lambda = 0},
\end{equation}
where $\deform{s}, \deform{s'} \in \deform{\mathcal{L}}$ are arbitrary
representatives of $s, s' \in \mathcal{L}$, respectively. It is well
known that this gives a well-defined full inner product for the
classical limit as we are in the unital situation and $\deform{h}$ was
non-degenerate and full already.
\begin{remark}
    \label{remark:LineBundleAlwaysHermitianMetric}%
    Of course, if $M$ is connected, then a pseudo Hermitian metric $h$
    on a \emph{line} bundle is necessarily either positive or negative
    definite. The issue of interesting signatures does not arise in
    this situation.
\end{remark}
Moreover, if $\deform{h}$ was completely positive then also its
classical limit is completely positive and hence a Hermitian fiber
metric, see the discussion in
\cite{bursztyn.waldmann:2005b}. Moreover, we note that the inner
product $\deform{h}$ as well as its classical limit $h =
\cl(\deform{h})$ do \emph{not} change when we move from
$\deform{\mathcal{L}}$ to ${}^\Psi\deform{\mathcal{L}}$: indeed, the
twisting with the symplectomorphism modifies only the \emph{left}
module structure. Thus we can absorb the symplectomorphism from the
beginning and assume $\deform{\mathcal{L}} = \Schnitte(L)[[\lambda]]$
also in this situation.  This induces an inner product $\deform{h}$ on
$\Schnitte(L)[[\lambda]]$ which can be written as
\begin{equation}
    \label{eq:DeformhAsSeries}
    \deform{h} = \sum_{r=0}^\infty \lambda^r h_r
\end{equation}
with sesquilinear maps $h_r: \Schnitte(L) \times \Schnitte(L)
\longrightarrow \Fkt$. In this case the classical limit of
$\deform{h}$ is, of course, given by $h_0$.

The following simple lemma clarifies the compatibility between the
classical limit of the inner product and the induced connection:
\begin{lemma}
    \label{lemma:connection in semiclassical limit}
    Let $\deform{\mathcal{L}} = (\Schnitte(L)[[\lambda]], \bullet,
    \bullet', \deform{h})$ be a $^\ast$-equivalence bimodule for the
    Hermitian star products $\star$ and $\star'$ deforming a pseudo
    Hermitian line bundle $(L, h_0)$. Assume furthermore $C_1 = C_1'$
    for the first order terms of the star products.  Then the induced
    connection $\nabla^L$ on $L$ is compatible with the induced pseudo
    Hermitian fiber metric $h_0$.
\end{lemma}
\begin{proof}
    The proof is a matter of straightforward computation.  Let $s, s'
    \in \Schnitte(L)$ and $f \in \Fkt$. The semiclassical limit of
    $\deform{h}(s, b \bullet f) = \deform{h}(s, b) \star f$ yields
    \[
    h_1(s, s'f) + h_0(s, R_1(s' ,f))
    =
    h_1(s, s') f + C_1(h_0(s, s'), f).
    \tag{$*$}
    \]
    From $\deform{h}(s \bullet f, s') = \cc{f} \star \deform{h}(s, s')$
    we get
    \[
    h_1(sf, s') + h_0(R_1(s,f), s')
    =
    \cc{f} h_1(s, s') + C_1(\cc{f},h_0(s, s')),
    \tag{$**$}
    \]
    and $\deform{h}(f \bullet' s, s') = \deform{h}(s, \cc{f} \bullet'
    s')$ gives
    \[
    h_1(fs, s') + h_0(R_1'(f, s), s')
    =
    h_1(s, \cc{f} s') + h_0(s, R_1'(\cc{f}, s')).
    \tag{$**$$*$}
    \]
    Taking $f = \cc{f}$ and substracting ($*$) from ($**$) yields
    \[
    h_1(sf, s') - h_1(s, fs')
    =
    h_0(s, R_1(s', f)) -h_0(R_1(s, f), s') + C_1(f, h_0(s, s')) -
    C_1(h_0(s, s'), f),
    \]
    whereas ($**$$*$)  shows, that 
    \[
    h_1(sf, s') - h_1(s, fs')
    =
    h_0(s, R_1'(f, s')) - h_0(R_1'(f, s), s').
    \]
    Together with $C_1(f, g) - C_1(g, f) = \I \{f, g\}$, we finally
    conclude that
    \[
    \Lie_{X_f}(h_0(s, s'))
    = \{h_0(s, s'), f\}
    = h_0\left(\nabla^L_{X_f}s, s'\right)
    + h_0\left(s, \nabla^L_{X_f}s'\right).
    \]
    Again, since the Hamiltonian vector fields span the tangent space
    at each point, this is enough to show the claim.
\end{proof}

After this preparation we want to discuss the semiclassical limit of
the $\lie{g}$-symmetry.  The following observation is now a trivial
computation:
\begin{lemma}
    \label{lemma:induced action}
    Let $\deform{\mathcal{L}}$ be a $\lie{g}$-equivariant equivalence
    bimodule. For $\xi \in \lie{g}$ and $s \in \mathcal{L} =
    \cl(\deform{\mathcal{L}})$,
    \begin{equation}
        \label{eq:classicalLimitAction}
        \rho(\xi)s
        = 
        \cl\left(\xi \acts \deform{s}\right)
    \end{equation}
    defines a lift of the classical $\lie{g}$-action on $\Fkt$ to
    $\mathcal{L}$ where $\deform{s} \in \deform{\mathcal{L}}$ is any
    representative of $s$.  This turns $\mathcal{L}$ into a
    $\mathcal{U}_{\mathbb{C}}(\lie{g})$-equivariant equivalence
    bimodule for $\Fkt$.
\end{lemma}
\begin{lemma}
    \label{lemma:PsiIsgInvariant}%
    Let $\deform{\mathcal{L}}$ be a $\lie{g}$-equivariant equivalence
    bimodule and let $\Psi$ be the unique symplectomorphism which
    makes the classical limit of ${}^\Psi\deform{\mathcal{L}}$ a
    static equivalence bimodule. Then $\Psi$ is $\lie{g}$-equivariant.
\end{lemma}
\begin{proof}
    We know that the classical limit $\cl(\deform{\mathcal{L}})$ of
    $\deform{\mathcal{L}}$ is a
    $\mathcal{U}_{\mathbb{C}}(\lie{g})$-equivariant equivalence
    bimodule for $\Fkt$. By Proposition~\ref{proposition:semidirect
      product for starPicH} we know that $\cl(\deform{\mathcal{L}})$
    splits uniquely into an equivariant automorphism $\Psi^*$ of
    $\Fkt$ and an equivariant static equivalence bimodule. This
    automorphism $\Psi^*$ is the pull-back of the symplectomorphism we
    are looking for.
\end{proof}
\begin{lemma}
    \label{lemma:InducedConnectionInvariant}%
    Let $\deform{\mathcal{L}} = (\Schnitte(L)[[\lambda]], \bullet,
    \bullet')$ be a $\lie{g}$-equivariant equivalence bimodule for the
    $\lie{g}$-invariant star products $\star$ and $\star'$ deforming a
    line bundle $L$. Assume $C_1 = C_1'$. Then the induced connection
    $\nabla^L$ is $\lie{g}$-invariant with respect to the classical
    limit action \eqref{eq:classicalLimitAction}.
\end{lemma}
\begin{proof}
    Let $\xi \in \lie{g}$, $f \in \Fkt$, and $s \in
    \Schnitte(L)$. Evaluating the equation $\xi \acts(s \bullet f) =
    (\xi \acts s) \bullet f + s \bullet (\xi \acts f)$ in order
    $\lambda$ gives $\rho(\xi) R_1(s, f) = R_1(\rho(\xi)s, f) + R_1(s,
    \rho(\xi)f)$.  Together with the analogous equation for the left
    module structure we get
    \[
    \rho(\xi)\left(\nabla^{L}_{X_{f}} s\right)
    =
    \nabla^L_{\rho(\xi)X_{f}} s + \nabla^L_{X_{f}} (\rho(\xi)s).
    \]
    As usual, this yields the desired compatibility for all vector
    fields.
\end{proof}
\begin{lemma}
    \label{lemma:InvariantHermitianMetric}%
    Let $\deform{\mathcal{L}} = (\Schnitte(L)[[\lambda]], \bullet,
    \bullet', \deform{h})$ be a $\lie{g}$-equivariant $^*$-equivalence
    bimodule for the $\lie{g}$-invariant Hermitian star products
    $\star$ and $\star'$ deforming a pseudo Hermitian line bundle $(L,
    h_0)$. Then the pseudo Hermitian fiber metric $h_0$ is
    $\lie{g}$-invariant with respect to the classical limit action
    \eqref{eq:classicalLimitAction}.
\end{lemma}
\begin{proof}
    First note that \eqref{eq:2kompatiblilaet H-wirkung} means $\xi
    \acts \deform{h}(s, s') = \deform{h}(\xi \acts s, s') +
    \deform{h}(s, \xi \acts s')$ for $\xi \in \lie{g}$ and $s, s' \in
    \Schnitte(L)$. Evaluating this in zeroth order gives immediately
    the result $\rho(\xi) h_0(s, s') = h_0(\rho(\xi) s, s') + h_0(s,
    \rho(\xi) s')$.  But this is the desired invariance of $h_0$.
\end{proof}
We can now collect the single results into the following theorem which
clarifies the semiclassical limit of $\lie{g}$-equivariant equivalence
bimodules completely:
\begin{theorem}[Semiclassical limit]
    \label{theorem:SemiclLimitWithAllStructures}%
    Assume that on $M$ we have a $\lie{g}$-invariant connection. Let
    $\star$ and $\star'$ be two $\lie{g}$-invariant star products on
    $M$.
    \begin{theoremlist}
    \item \label{item:RingVersionSCL} Assume $\star$ and $\star'$ are
        $\lie{g}$-equivariantly Morita equivalent via a
        $\lie{g}$-equivariant equivalence bimodule
        $\deform{\mathcal{L}}$. Then there exists a unique
        $\lie{g}$-equivariant symplectomorphism $\Psi$, a line bundle
        $L$ unique up to isomorphism with a lift of the
        $\lie{g}$-action, and a $\lie{g}$-invariant connection
        $\nabla^L$ on $L$ such that $\deform{\mathcal{L}}$ is
        isomorphic to the $\Psi$-twist of a bimodule deformation
        $(\Schnitte(L)[[\lambda]], \bullet, \bullet'')$ of the
        sections of the line bundle with respect to $\star$ and
        $\star''$ where $\star''$ is $\lie{g}$-equivariantly
        equivalent to $\Psi^*(\star')$ such that $\star''$ has the
        same first order term than $\star$. Finally, the connection is
        obtained as first order difference of the deformed bimodule
        deformations $\bullet$ and $\bullet''$.
    \item \label{item:StarVersionSCL} If in addition $\star$ and
        $\star'$ are Hermitian and $\deform{\mathcal{L}}$ is even a
        $\lie{g}$-equivariant $^*$-equivalence bimodule then in
        addition we have a $\lie{g}$-invariant pseudo Hermitian metric
        $h_0$ on $L$, with $\nabla^L$ being metric, such that the
        $(\Fkt[[\lambda]], \star)$-valued inner product is a
        deformation of $h_0$.
    \item \label{item:StrongVersionSCL} If furthermore
        $\deform{\mathcal{L}}$ is even a $\lie{g}$-equivariant strong
        equivalence bimodule then $h_0$ is Hermitian.
    \end{theoremlist}
\end{theorem}
\begin{proof}
    Let $\deform{\mathcal{L}}$ be the $\lie{g}$-equivariant
    equivalence bimodule for $\star'$ acting from the left and $\star$
    from the right. First we apply
    Proposition~\ref{proposition:clLimitEquivalenceBimoduL} to split
    off the symplectomorphism $\Psi$ by passing from $\star'$ to
    $\Psi^*(\star')$. By Lemma~\ref{lemma:PsiIsgInvariant}, $\Psi$ is
    necessarily $\lie{g}$-equivariant so $\Psi^*(\star')$ is still a
    $\lie{g}$-invariant star product. The twisted bimodule is now a
    $\lie{g}$-equivariant equivalence bimodule for $\Psi^*(\star')$
    and $\star$ such that its classical limit is \emph{static}. Next
    we can pass from $\Psi^*(\star')$ to an equivalent star product
    via a $\lie{g}$-equivariant equivalence transformation in such a
    way that the new star product $\star''$ has the same first order
    term as $\star$.  Here we need the $\lie{g}$-invariant connection
    on $M$ to achieve this in a $\lie{g}$-equivariant way: indeed,
    having a $\lie{g}$-invariant connection, we can write the
    $\lie{g}$-invariant Hochschild $2$-coboundary $C_1 - \Psi^*(C_1')$
    as Hochschild $\delta$ of some \emph{$\lie{g}$-invariant}
    differential operator $D_1$ according to
    \cite[Prop.~2.1]{bertelson.bieliavsky.gutt:1998a}. Then
    exponentiating $D_1$ gives a $\lie{g}$-equivariant equivalence
    transformation from $\Psi^*(\star')$ to $\star''$ where now
    $\star''$ has the same first order as $\star$.  Correspondingly,
    we can transport the equivalence bimodule along. This results in a
    $(\star'', \star)$-equivalence bimodule which is isomorphic to
    $\Schnitte(L)[[\lambda]]$ endowed with deformed bimodule
    multiplications $\bullet$ and $\bullet''$ where $L$ is a line
    bundle uniquely determined up to isomorphism. Moreover, by
    Lemma~\ref{lemma:induced action} the line bundle carries an action
    of $\lie{g}$ which is a lift of the action on $M$. With
    Lemma~\ref{lemma:connection in semiclassical limit} we get the
    connection $\nabla^L$ which is $\lie{g}$-invariant by
    Lemma~\ref{lemma:InducedConnectionInvariant}. This shows the
    ring-theoretic part. For the case of Hermitian star products and
    $^*$-equivalences we use Lemma~\ref{lemma:connection in
      semiclassical limit} and
    Lemma~\ref{lemma:InvariantHermitianMetric} to conclude that the
    induced pseudo Hermitian metric has the desired
    properties. Finally, the last part is clear.
\end{proof}
\begin{remark}
    \label{remark:NoNeedForConnections}%
    If we are only interested in star products where the first order
    term is $\frac{\I}{2} \{\cdot, \cdot\}$ (and not just its
    antisymmetric part) then the extra assumption of an invariant
    connection on $M$ becomes superfluous. Also, if we are not
    interested in the invariant connection $\nabla^L$ on $L$ we can
    omit this assumption.
\end{remark}

%
%

\subsection{Characterization of equivariant Morita equivalent star
  products}
\label{subsec:EquivariantMEStarProducts}%

Now we look for necessary and sufficient conditions for equivariant
Morita equivalence of star products: in some sense we want to
investigate how the conditions in the semiclassical limit are already
sufficient to guarantee the corresponding equivalence also in the full
``quantum'' situation.

It is well known that two star products $\star$ and $\star'$ over $(M,
\omega)$ are Morita equivalent if and only if there is a line bundle
$L \longrightarrow M$ and a symplectomorphism $\Psi: M \longrightarrow
M$ such that $\Psi^{\ast} c(\star')- c(\star) = 2\pi \I c_1(L)$, see
\cite{bursztyn.waldmann:2002a}.  As before, $\Psi^{\ast}(\star')$
denotes the star product $\star'$ twisted with $\Psi$ and $c(\star)$
is the characteristic class of $\star$.  Finally, $c_1(L) \in
\HdR^2(M, \mathbb{Z})$ is the \emph{Chern class} of $L$ viewed as
class in the deRham cohomology.  Choosing a covariant derivative
$\nabla^L$ with curvature two-from $R^L \in \Formen^2(M)$ for $L$, we
have $[R^L] = 2 \pi \I c_1(L) \in \HdR^2(M, \mathbb{C})$. Note that
this is \emph{not} the Chern class in the second integral cohomology
of $M$ but it's image in the deRham cohomology, so possible torsion
bundles will have trivial class $c_1(L) = 0$ in this version.

First we want to assume that there is already an action of $\lie{g}$
on a given line bundle $L$. Then we can modify the Fedosov
construction for the bimodule $\Schnitte(L)$ to the
$\lie{g}$-invariant case, as is done in \cite{jansen:2006a,
  waldmann:2002b} for the case of a Lie group action on $L
\longrightarrow M$.
\begin{proposition}[Invariant Fedosov construction]
    \label{proposition:Fedosov for equiva morita}%
    Let $L \longrightarrow M$ be a line bundle over a symplectic
    manifold $M$. Suppose $\lie{g}$ acts on $M$ via symplectic vector
    fields and there is a lift of this action to $L$. Moreover,
    suppose there is a $\lie{g}$-invariant symplectic connection on
    $M$, $\Omega \in \lambda \Formen^2(M, \mathbb{C})[[\lambda]]$ is a
    formal series of $\lie{g}$-invariant closed two-forms, and $L$
    allows for a $\lie{g}$-invariant connection $\nabla^L$.
    \begin{propositionlist}
    \item \label{item:InvariantFedosov} The Fedosov construction leads
        to a $\lie{g}$-equivariant equivalence bimodule
        $(\Schnitte(L)[[\lambda]], \bullet, \bullet')$ for the
        $\lie{g}$-invariant Fedosov star products $\star_{\nabla,
          \Omega}$ and $\star' = \star_{\nabla, \Omega'}$ with $\Omega'
          = \Omega + \I\lambda R^L$, acting from the left. Hence the
          $\lie{g}$-invariant characteristic classes are
          $c^{\lie{g}}(\star_{\nabla, \Omega}) = \frac{1}{\I
            \lambda}[\omega + \Omega]$ and $c^{\lie{g}}(\star') =
          c^{\lie{g}}(\star) + 2 \pi \I c_1(L)$.
      \item \label{item:InvariantStarFedosov} If in addition $\Omega$
          is real and if there exists a $\lie{g}$-invariant pseudo
          Hermitian metric $h_0$ on $L$ such that $\nabla^L$ is metric
          then $\star$ and $\star'$ are Hermitian and $h_0$ can be
          deformed into a $\lie{g}$-invariant $(\Fkt[[\lambda]],
          \star)$-valued inner product making the above bimodule even
          a $\lie{g}$-equivariant $^*$-equivalence bimodule.
    \item \label{item:InvariantStrongFedosov} In in addition $h_0$ is
        Hermitian then the resulting inner products are completely
        positive and the bimodule is a $\lie{g}$-equivariant strong
        equivalence bimodule.
    \end{propositionlist}
\end{proposition}
\begin{proof}
    The functoriality of the Fedosov construction \cite{jansen:2006a,
      waldmann:2002b} give this result immediately.
\end{proof}
\begin{lemma}
    \label{lemma:liftingwithCartan-Model}
    Let $L \longrightarrow M$ be a line bundle such that $c_1(L) \in
    \HdR^2(M, \mathbb{Z})$ is in the image of the combined map of
    \eqref{eq:EquivariantToInvariantdRh} after
    \eqref{eq:CanonicalHdRgToHdR}. For every Hermitian fiber
    metric $h_0$ on $L$ there is a metric connection $\nabla^L$ and a
    lifting of the $\lie{g}$-action to $L$ such that $h_0$ and
    $\nabla^L$ are invariant.
\end{lemma}
From this folklore statement, see e.g. \cite{kostant:1970a,
  lashof:1991a}, we immediately obtain the characterization of
equivariant (strong) Morita equivalence proving
Theorem~\ref{theorem:TheRealStuff}:
\begin{proof}[ of Theorem~\ref{theorem:TheRealStuff}]
    If they are $\lie{g}$-equivariantly (strongly) Morita equivalent
    then our considerations from the semiclassical limit in
    Theorem~\ref{theorem:SemiclLimitWithAllStructures} show that their
    characteristic classes combine into the $\lie{g}$-invariant Chern
    class of a line bundle $L$ with $\lie{g}$-action, $2 \pi \I c_1(L)
    = \Psi^*c(\star') - c(\star) \in 2 \pi \I \HdR(M,
    \mathbb{Z})$. Without $\lie{g}$-invariance, this is the original
    characterization from \cite{bursztyn.waldmann:2002a}. Moreover, we
    know that on $L$ there is a $\lie{g}$-invariant connection
    $\nabla^L$ and hence its curvature $R^L$ is a $\lie{g}$-invariant
    representative of $2 \pi \I c_1(L)$. Now define $J^L(\xi) \in
    \End_{C^\infty(M)}(\Schnitte(L)) = C^\infty(L)$ via $J^L(\xi) = -
    \rho(\xi) - \nabla^L_{X_\xi}$ where $\rho(\xi)$ is the action of
    $\xi \in \lie{g}$ on $\Schnitte(L)$. From the definition of
    $\rho(\xi) = - \Lie_{X_\xi}$ on functions it follows that
    $J^L(\xi)$ is indeed $C^\infty(M)$-linear and hence a function
    itself. Since $\nabla^L$ is invariant and since $\rho(\cdot)$ is a
    Lie algebra action by module derivations it follows by a
    straightforward computation that $J^L$ extends $R^L$ to an
    equivariant two-form as wanted.

    Conversely, we first note that every $\lie{g}$-invariant
    (Hermitian) star product is $\lie{g}$-equivariantly equivalent
    (with a \emph{real} equivalence transformation since a
    $\lie{g}$-equivariant equivalence can be made a $^*$-equivalence
    according to \cite[Exercise~7.6]{waldmann:2007a} without
    destroying the $\lie{g}$-equivariance) to a Fedosov star product
    (with real $\Omega$).  The first condition on the characteristic
    classes shows that the star products are Morita equivalent and the
    corresponding line bundle can be used to build the equivalence
    bimodule. Note that $\Psi^*c(\star') - c(\star)$ is necessarily
    concentrated in zeroth order of $\lambda$ only.  By
    Lemma~\ref{lemma:liftingwithCartan-Model} we have a lift to this
    line bundle, a $\lie{g}$-invariant Hermitian metric, and a
    $\lie{g}$-invariant metric connection. Then the Fedosov
    construction from Proposition~\ref{proposition:Fedosov for equiva
      morita} ensures that the corresponding Fedosov star products and
    hence also the original star products are $\lie{g}$-equivariantly
    (strongly) Morita equivalent.
\end{proof}

%
%

\subsection{The role of quantum momentum maps}
\label{subsec:RoleOfQuantumMomentumMaps}

We conclude this section with some remarks on the role of quantum
momentum maps.

To apply our lifting statements for general Hopf algebra actions from
Theorem~\ref{theorem:splitting the image starpich to starpic} we have
to use the Hopf $^*$-algebra $\mathcal{U}_\lambda(\lie{g})$ for star
products with quantum momentum map instead of
$\mathcal{U}_{\mathbb{C}}(\lie{g})$. In this case, the action on
$C^\infty(M)[[\lambda]]$ of $\mathcal{U}_\lambda(\lie{g})$ is given by
\eqref{eq:xiactsf} and therefore starts in \emph{first} order of
$\lambda$ but not in the zeroth order.  We can always pass from a
$\mathcal{U}_{\mathbb{C}}(\lie{g})$-action to a corresponding action
of $\mathcal{U}_\lambda(\lie{g})$ in the following way: Let
$\deform{\mathcal{L}}$ be a $\mathbb{C}[[\lambda]]$-module with a left
$\mathcal{U}_{\mathbb{C}}(\lie{g})$-action $\rho$ on it. Then 
\begin{equation}
    \label{eq:FromUtoUlambda}
    \xi \acts x = \I \lambda \rho(\xi) x
\end{equation}
for $\xi \in \lie{g}$ and $x \in \deform{\mathcal{L}}$ extends
uniquely to an action of $\mathcal{U}_\lambda(\lie{g})$ on
$\deform{\mathcal{L}}$.  However, not every
$\mathcal{U}_\lambda(\lie{g})$-action is of that form. The reason is
that \emph{all} operators $\xi \acts \cdot$ have to vanish for
$\lambda = 0$, i.e. in the classical limit, for $\xi \in \lie{g}$.
However, for the lift according to Theorem~\ref{theorem:splitting the
  image starpich to starpic} to equivalence bimodules this is
automatically the case:
\begin{lemma}
    \label{lemma:WhichUlambdaActionsAreCool}%
    Let two $\lie{g}$-invariant star products $\star$, $\star'$ be
    Morita equivalent with an equivalence bimodule
    $\deform{\mathcal{L}}$, and let $\Psi$ be the symplectomorphism
    such that ${}^\Psi \cl(\deform{\mathcal{L}})$ is static. Assume
    $\star$ and $\star'$ have quantum momentum maps $J$ and $J'$.  For
    the canonical lift \eqref{eq:LiftedAction} of the
    $\mathcal{U}_\lambda(\lie{g})$-action to $\deform{\mathcal{L}}$
    there exists a unique action $\rho$ of
    $\mathcal{U}_{\mathbb{C}}(\lie{g})$ such that $\acts$ is of the
    form \eqref{eq:FromUtoUlambda} iff the induced symplectomorphism
    $\Psi$ satisfies $\Psi^* J'_0 = J_0$.
\end{lemma}
\begin{proof}
    We know that $\deform{\mathcal{L}} \cong
    {}^\Psi\Schnitte(L)[[\lambda]]$ as $\mathbb{C}[[\lambda]]$-modules
    which induces deformed module structures $\bullet'$ and $\bullet$
    on ${}^\Psi\Schnitte(L)[[\lambda]]$. The lifted action
    \eqref{eq:LiftedAction} is now given by
    \[
    \xi \acts s = J'(\xi) \bullet' s - s \bullet J(\xi)
    \tag{$*$}
    \]
    for $s \in \Schnitte(L)[[\lambda]]$ and $\xi \in \lie{g}$. The
    classical limit of ($*$) is $\Psi^*(J'_0(\xi)) s - s J_0(\xi)$
    which is zero iff $\Psi^*J'_0(\xi) = J_0(\xi)$ for all $\xi \in
    \lie{g}$. It is clear that ($*$) is of the form
    \eqref{eq:FromUtoUlambda} iff it vanishes in the zeroth order.
\end{proof}

Since both classical limits $J_0$ and $J'_0$ are classical momentum
maps it follows that a symplectomorphism $\Psi$ with $\Psi^* J'_0 =
J_0$ is $\lie{g}$-equivariant. Conversely, if $\Psi$ is
$\lie{g}$-equivariant then $\Psi^* J'_0$ is still a classical momentum
map and hence differs from $J_0$ only by a cocycle $c$ in
$Z^1(\lie{g}, \mathbb{C})$. Then the lifted
$\mathcal{U}_\lambda(\lie{g})$-action in zeroth order is given by $\xi
\acts s = c_\xi s + \cdots$ which does \emph{not} vanish.  However, in
this case we can modify the $\mathcal{U}_\lambda(\lie{g})$-action by
precisely this cocycle. This results in a
$\mathcal{U}_\lambda(\lie{g})$-equivalence bimodule with action
\begin{equation}
    \label{eq:FunnyNewAction}
    \xi \acts' s = J'(\xi) \bullet' s - s \bullet J(\xi) - c_\xi s
\end{equation}
for $\xi \in \lie{g}$ and $s \in \Schnitte(L)[[\lambda]]$. By
$c_{[\xi, \eta]} = 0$ this clearly extends to a
$\mathcal{U}_\lambda(\lie{g})$-action. Since $c_\xi$ is a
\emph{constant} we still have compatibility with the bimodule
structure. This way, the action $\acts'$ is now of the form
\eqref{eq:FromUtoUlambda} with an action $\rho$ of
$\mathcal{U}_{\mathbb{C}}(\lie{g})$. We summarize this in the
following lemma:
\begin{lemma}
    \label{lemma:UlambdaUEME}%
    Let $\star$ and $\star'$ be $\lie{g}$-invariant star products with
    quantum momentum maps $J$ and $J'$. Let $\deform{\mathcal{L}}$ be
    an equivalence bimodule with induced symplectomorphism $\Psi$. If
    $\Psi$ is $\lie{g}$-equivariant then $\star$ and $\star'$ are
    $\lie{g}$-equivariantly Morita equivalent.
\end{lemma}
\begin{proof}
    Indeed, by the above considerations the cocycle $c_\xi =
    \Psi^*J'_0(\xi) - J_0(\xi)$ can be used to modify the lifted
    action \eqref{eq:LiftedAction} to obtain an action of the form
    \eqref{eq:FromUtoUlambda}. Thus $\deform{\mathcal{L}}$ with this
    action is a $\mathcal{U}_{\mathbb{C}}(\lie{g})$-equivariant
    equivalence bimodule for $\star'$ and $\star$.
\end{proof}

In this situation where we have such an equivalence bimodule with
$\lie{g}$-invariant symplectomorphism $\Psi$ we see that also
$\Psi^*(\star')$ is a $\lie{g}$-invariant star product with quantum
momentum map $\Psi^*J'$. Moreover, we see that
$c^{\lie{g}}(\Psi^*(\star')) = \Psi^* c^{\lie{g}}(\star')$ since
$\Psi$ is $\lie{g}$-equivariant. By Morita equivalence we have that
$\Psi^* c^{\lie{g}} (\star') - c^{\lie{g}}(\star)$ maps to the
$2\pi\I$-integral deRham cohomology, i.e. $\Psi^* c^{\lie{g}}(\star') - c^{\lie{g}}(\star)
\in 2 \pi \I\HdR^2(M, \mathbb{Z})$. Finally, since the star products
$\Psi^*(\star')$ and $\star$ have quantum momentum maps it is clear
that their invariant characteristic classes are in the image of
\eqref{eq:EquivariantToInvariantdRh} by
Proposition~\ref{proposition:existence of momentum maps and
  equivariant cohomology}. But then \emph{trivially} the difference
$\Psi^* c^{\lie{g}} (\star') - c^{\lie{g}}(\star)$ is in the image as
well. This way, we recover the general statement of
Theorem~\ref{theorem:TheRealStuff} for star products with quantum
momentum maps.  Note that now the lifting part of
Theorem~\ref{theorem:TheRealStuff} was now done using the general
lifting construction from Theorem~\ref{theorem:splitting the image
  starpich to starpic}.
\begin{remark}
    \label{remark:MOMAStar}%
    In the case of strong Morita equivalence we can argue analogously,
    observing that for real momentum maps the above cocycle is real,
    too. Hence the inner products will not be influenced at all.
\end{remark}

\begin{footnotesize}
\renewcommand{\arraystretch}{0.5} 

\end{footnotesize}

\end{document}